\documentclass[11pt,reqno,a4paper]{article}
\usepackage{geometry}
\usepackage[blocks]{authblk}
\usepackage{amsmath,amsthm,amstext,amscd,amssymb,euscript,url,mathrsfs,euscript,dsfont,calrsfs,mathtools,comment,relsize,xcolor}
\usepackage{hyperref}
\usepackage{amsmath,amsthm,amstext,amscd,amssymb,euscript,url}
\usepackage{mathrsfs}
\usepackage{calrsfs}
\usepackage{epsfig}
\usepackage[inline]{enumitem}
\usepackage{todonotes}
\usepackage{mathtools}
\usepackage[absolute,overlay]{textpos}
\usepackage{graphicx}
\usepackage{subfig}
\usepackage{tikz}
\usetikzlibrary{matrix}
\usetikzlibrary{arrows,positioning}

\usepackage{color}

\newcommand{\R}{\mathbb R}
\newcommand{\N}{{\mathbb N}}

\newcommand{\E}{\mathbb E}

\renewcommand{\phi}{\varphi}

\newcommand{\pee}{\ensuremath{\mathbb{P}}}

\newcommand{\loc}{\mathcal{L}}

\def\1{{\mathchoice {\rm 1\mskip-4mu l} {\rm 1\mskip-4mu l}
{\rm 1\mskip-4.5mu l} {\rm 1\mskip-5mu l}}}

\newtheorem{theorem}{{\small T}{\scriptsize HEOREM}}[section]
\newtheorem{corollary}{{\bf{\small C}{\scriptsize OROLLARY}}}[section]
\newtheorem{proposition}{{\bf{\small P}{\scriptsize ROPOSITION}}}[section]
\newtheorem{lemma}{{\bf{\small L}{\scriptsize EMMA}}}[section]
\newtheorem{remark}{{\bf{\small R}{\scriptsize EMARK}}}[section]
\newtheorem{definition}{{\bf{\small D}{\scriptsize EFINITION}}}[section]
\newtheorem{induction}{{\bf{\small I}{\scriptsize NDUCTIVE HYPOTHESIS}}}[section]

\renewenvironment{proof}[1]
{\noindent{{\bf{\small{ P}{\scriptsize ROOF}}}.}\hspace{0.1cm} #1} {$\;\qed$\newline}

\newcommand{\beq}{\begin{eqnarray}}
\newcommand{\eeq}{\end{eqnarray}}

\newcommand{\ba}{\begin{align*}}
\newcommand{\ea}{\end{align*}}

\newcommand{\be}{\begin{equation}}
\newcommand{\ee}{\end{equation}}

\newcommand{\bl}{\begin{lemma}}
\newcommand{\el}{\end{lemma}}

\newcommand{\br}{\begin{remark}}
\newcommand{\er}{\end{remark}}

\newcommand{\bt}{\begin{theorem}}
\newcommand{\et}{\end{theorem}}

\newcommand{\bd}{\begin{definition}}
\newcommand{\ed}{\end{definition}}

\newcommand{\bind}{\begin{induction}}
\newcommand{\eind}{\end{induction}}

\newcommand{\bp}{\begin{proposition}}
\newcommand{\ep}{\end{proposition}}

\newcommand{\bc}{\begin{corollary}}
\newcommand{\ec}{\end{corollary}}

\newcommand{\bpr}{\begin{proof}}
\newcommand{\epr}{\end{proof}}

\newcommand{\bi}{\begin{itemize}}
\newcommand{\ei}{\end{itemize}}

\newcommand{\ben}{\begin{enumerate}}
\newcommand{\een}{\end{enumerate}}

\newcommand{\caB}{{\mathcal B}}
\newcommand{\caC}{{\mathscr C}}

\newcommand{\caD}{{\EuScript D}}

\newcommand{\caG}{{\mathcal G}}

\newcommand{\caH}{{\mathcal H}}

\newcommand{\caL}{{\mathcal L}}
\newcommand{\caM}{{\mathcal M}}

\newcommand{\caS}{{\mathcal S}}

\newcommand{\Beta}{\text{Beta}}

\newmuskip\pFqmuskip

\newcommand\pFq[6][8]{%
	\begingroup 
	\pFqmuskip=#1mu\relax
	\mathcode`\,=\string"8000
	\begingroup\lccode`\~=`\,
	\lowercase{\endgroup\let~}\pFqcomma
	{}_{#2}F_{#3}{\left[\genfrac..{0pt}{}{#4}{#5};#6\right]}%
	\endgroup
}
\newcommand{\pFqcomma}{\mskip\pFqmuskip}

\newtheorem{asu}{Assumption}
\newcommand{\basu}{\begin{asu}}
\newcommand{\easu}{\end{asu}}

\begin{document}
\title{Intertwining and propagation of mixtures for\\
generalized KMP models and harmonic models}
\author[1]{Cristian Giardin\`a
\thanks{Email: \texttt{cristian.giardina@unimore.it}}
}
\author[2]{Frank Redig
\thanks{Email: \texttt{F.H.J.Redig@tudelft.nl}}
}
\author[2]{Berend van Tol
\thanks{Email: \texttt{B.T.vanTol@tudelft.nl}}
}

\affil[1]{ FIM, University of Modena and Reggio Emilia, Modena, Italy}
\affil[2]{Institute of Applied Mathematics, Delft University of Technology, Delft,
The Netherlands}

\date{Dated: \today}

\maketitle

\begin{abstract}
We study a class of stochastic models of  mass transport on discrete vertex set $V$. For these models, a one-parameter family of homogeneous product measures
$\otimes_{i\in V} \nu_\theta$ is reversible.
We prove that the set of mixtures of inhomogeneous product measures with equilibrium marginals, i.e., the set of measures of the form
\[
\int\Big(\bigotimes_{i\in V} \nu_{\theta_i}\Big) \,\Xi(\prod_{i\in V}d\theta_i)
\]
is left invariant by the dynamics in the course of time, and the  ``mixing measure'' $\Xi$ evolves according to a Markov process which we then call ``the hidden parameter model''.
This generalizes results from \cite{deMasi} to a larger class of models and on more general graphs. The class of models includes discrete and continuous generalized KMP models, as well as discrete and continuous harmonic models.
The results imply that in all these models, the non-equilibrium steady state of their reservoir  driven version is a mixture of product measures where the mixing measure is in turn the stationary state of the corresponding ``hidden parameter model''. For the boundary-driven harmonic models on the chain
$\{1,\ldots, N\}$ with nearest neighbor edges, we recover that the stationary measure of the hidden parameter model is the joint distribution of the ordered Dirichlet distribution (cf.\ \cite{red}), 
with a purely probabilistic proof  based on a spatial Markov property of the hidden parameter model.
\end{abstract}

\newpage

\section{Introduction}
Recent developments in the study of the KMP model 
and related models have revealed that the non-equilibrium steady state of the boundary driven version of such models is a mixture of product measures of equilibrium marginals. In the simplest setting of the KMP model \cite{deMasi}, this means that the non-equilibrium steady state is a mixture of products of exponential distributions, where the joint distribution of the parameters of these exponentials is in turn a stationary distribution of an auxiliary model, the so-called {\em hidden temperature model} \cite{deMasi}.
For a related  class of models, the generalized  harmonic models \cite{frassek2020non}, \cite{frass1}, \cite{frass2},   the non-equilibrium steady state of the continuous model is given in closed form in terms of products of gamma distributions, with identical shape parameters, and where the scale parameters have the ordered Dirichlet distribution \cite{red}. In the simplest setting of the harmonic model, the non-equilibrium steady state is a product of exponential distributions, where the (scale) parameters are distributed as the order statistics of i.i.d. uniforms \cite{gab}. The structure of the stationary state as a mixture was already conjectured in \cite{bert} (for the KMP model),  based on  macroscopic fluctuation theory.

So far, these results are all obtained in the setting of a chain geometry, with boundary reservoirs at left and right ends. Most of these results are strongly based on dualities, which reduce the computation of moment of order $n$ in the non-equilibrium steady state to the computation of absorption probabilities of $n$ dual particles. For the characterization of the non-equilibrium steady state of the generalized harmonic models of a chain, an additional input came from integrability. Is it usually the latter which provides closed-form expressions for the absorption probabilities of the dual process and is only applicable in the chain geometry, whereas duality results are valid in a setting of general graphs.

In this paper, using a reformulation of duality as an {\em intertwining} relation, we prove that for a large class of models {\em on a general graph}, there exist {\em hidden parameter models}. As a consequence, the non-equilibrium steady state is a mixture of equilibrium product marginals where the mixing measure (i.e., the joint distribution of the parameters of these marginals) is the stationary measure of the corresponding hidden parameter model. This stationary measure is usually inaccessible in explicit form on a general graph. In the case of the harmonic model on a chain, we are able to prove that it coincides with the mixing measure found in \cite{red}, using probabilistic arguments only based on a Markovian structure in the hidden parameter model.
These results show that in essence, the existence of hidden parameter models is based on duality, and therefore not restricted to integrable models.
However, the identification of the mixing measure, i.e., the measure describing the joint distribution of the parameters, is only possible when there is extra structure (i.e., extra symmetries) which makes it possible e.g. to use the quantum inverse scattering method
\cite{frass1}, or put more probabilistically, to have a Markovian structure of the mixing measure.

The rest of our paper is organized as follows.
In section \ref{genst}, we sketch the general structure of the models under consideration.
In section \ref{section2} we discuss the discrete and continuous generalized KMP models, recovering and generalizing the hidden temperature models in \cite{deMasi}. In section  \ref{section3} we deal with the generalized harmonic models. In particular, we identify the corresponding hidden parameter models
and establish intertwining on a generic graph.
For the boundary driven chain we characterize
the stationary measure of the hidden parameter model
by using a self-contained argument which rely on
the particular structure of the model.
In section  \ref{section4} we extend the analysis to another model, the symmetric inclusion process (SIP) and prove that it admits Poisson intertwining. As a consequence the non-equilibrium steady state is a mixture of Poisson product measures, where the mixture measure is a non-equilibrium steady state of a corresponding continuous model (the Brownian energy process). We also recover the simplest setting of boundary driven independent random walks, where the intertwined dynamics is deterministic and has a unique fixed point, which implies that the non-equilibrium steady state is a product of Poisson measures. The latter is of course well-known but we believe it is still insightful to recover it from the point of view of intertwining.
\subsection{General structure of the models}\label{genst}
We consider a finite set of vertices $V$, and a symmetric irreducible collection of edge weights
$p(i,j)=p(j,i)\geq 0$ where $i,j\in V$. Here, by irreducibility we mean that for every $i,j\in V$ there
exists a finite discrete path $\gamma(0), \ldots, \gamma (n)$ with $\gamma(0)=i, \gamma(n)=j$ and
$p(\gamma(i), \gamma(i+1))>0$ for all $i=0,\ldots, n-1$.
We will then consider Markov processes on either the state space $\N^V=\{0,1,2,\ldots\}^V$ (discrete models) or
the state space $\R_+^V=[0,\infty)^V$ (continuous models).
The generator of these processes will  take the form
\[
\sum_{i,j\in V} p(i,j) L_{ij}
\]
where $L_{ij}$ is the so-called single edge generator, which acts only on the variables $\eta_i,\eta_j$ and models the transport of mass  along the edge connecting the sites $i,j\in V$.
For the boundary driven version of the models we have a generator
taking the form
\[
\sum_{i,j\in V} p(i,j) L_{ij} + \sum_{i\in V} c(i)  L_{\theta_i^*}
\]
where $c(i) \ge 0$ is a non-negative constant tuning the coupling of site $i\in V$ to a ``reservoir'' with parameter $\theta_i^*>0$.
The single-site generator $L_{\theta_i^*}$ is acting only on the variables $\eta_i$ 
and models the input and output of mass at the vertex $i\in V$,
by fixing the average number of particles to $\theta_i^*$. 

The system with generator $\sum_{i,j\in V} p(i,j) L_{ij}$ will have a one parameter family of
product invariant measures $\bigotimes_{i\in V}\nu_\theta$, where the parameter $\theta>0$ labels the expected number
of particles (or mass) and corresponds to the conserved quantity (total number of particles or total mass).
Then the system coupled to reservoirs with identical parameters ($\theta_i^*=\theta^*$ for $i\in V$) has a unique stationary measure $\bigotimes_{i\in V}\nu_{\theta^*}$. If the reservoir parameters are different, then
the unique stationary measure is no longer a product measure, and is called a {\em non-equilibrium steady state}, where non-equilibrium refers to the absence of reversibility.

The main aim of this paper is to understand for a family of models of this type the propagation of {\em inhomogeneous} product measures 
$\bigotimes_{i\in V} \nu_{\theta_i}$ in the course of time.
Given $\theta=(\theta_i)_{i\in V}$,
we will then find that these measures are mapped to a stochastic mixture of the form
\[
\E_\theta \Big(\bigotimes_{i\in V} \nu_{\theta_i(t)}\Big)
\]
where $(\theta_i(t), t\geq 0, i\in V)$ will evolve as a Markov process which we then call, following \cite{deMasi}, the ``hidden parameter model''. As a consequence, the unique stationary measure (non-equilibrium steady state) will be a mixture of product measures of the type
\[
\int\Big(\bigotimes_{i\in V} \nu_{\theta_i}\Big) \Xi(\prod_{i\in V}d\theta_i)
\]
The ``mixing measure'' $\Xi$ is then the unique invariant measure of the hidden parameter model.
Thus, in the reservoir-driven setup, the identification of the non-equilibrium steady state is reduced to the
identification of the stationary measure of the hidden parameter model.

The two most important examples of models having the property that the set of mixture of equilibrium product measure is closed
under the dynamics  will be models of ``KMP type'' (section \ref{section2}) or models of ``harmonic type'' (section \ref{section3}).
For another class of models, namely the symmetric inclusion process and the independent random walkers (section \ref{section4}),
we will show that the same happens with a product of Poisson measures, where
the evolution of the Poisson parameters is then either a Markov diffusion process or a deterministic process.

In what follows we will always use an upright $L$ for the generator of the process under study,
and the symbol $\loc$ for the corresponding hidden parameter model.
We will always use the notation $\E_\eta, \E_\xi$ for expectations for process with discrete state space such as $\N^V$, $\E_\zeta$ for the expectations of processes with continuous state space such as $[0,\infty)^V$ and $\E_\theta$ for expectations of processes of hidden parameter models
(also with state space $[0,\infty)^V$).

\section{Generalized KMP processes}
\label{section2}

In this section we study the discrete, resp. continuous, generalized KMP models, parametrized by a non-negative number $s>0$. These models are a one-parameter generalization of the original KMP model and
were introduced in \cite{GKRV}.
For arbitrary $s>0$ we prove new dualities with a generalized hidden parameter model.
This in turn implies that products of discrete gamma, resp. continuous gamma,
distributions evolve in the course of time into mixtures of such product measures,
where the mixing measure is the distribution of the corresponding hidden parameter model.

We start by first considering the bulk process and then we add reservoirs.
The original discrete and continuous KMP models \cite{kmp} will be recovered for $s=1/2$.

\subsection{Discrete generalized KMP}

We consider a finite set of vertices $V$, and irreducible edge rates $p(i,j)$, as outlined in section \ref{genst}.
The discrete generalized KMP process with parameter $2s>0$ is a Markov process on $\N^V$ and is defined via the generator
\be\label{diskmp}
L f(\eta)
= \sum_{i,j\in V} p(i,j) L_{ij} f(\eta)\,.
\ee
Here the single edge generator $L_{ij}$ acts on the variables $\eta_i,\eta_j$ as 
\be\label{diskmpsedge}
L_{ij} f(x,y)= \E \big (f(X, x+y-X) - f(x,y)\big)
\ee
where $X$ is beta-binomial with parameters $x+y, 2s, 2s$, i.e.,
\be\label{betabin}
\pee(X=k)= \int_0^1 {x+y\choose k} p^k (1-p)^{x+y-k} \Beta(2s,2s) [dp]
\ee
where $k\in \{0,1,\ldots,x+y\}$ and
\be\label{beta}
\Beta(2s,2s) [dp]= \frac{1}{B(2s,2s)} p^{2s-1} (1-p)^{2s-1} dp
\ee
denotes the Beta distribution with parameters $(2s,2s)$.

The discrete generalized KMP process has reversible product measures which are product of discrete Gamma distributions parametrized as follows
\be\label{discgam}
\nu_\theta(n)= \left(\frac{\theta}{1+\theta}\right)^{n} \frac{\Gamma(2s+n)}{\Gamma(2s)}  \left(\frac{1}{1+\theta}\right)^{2s}
\ee
The relation between the parameter $\theta$ and the expectation of the marginals is given by
\be
\sum_{n=0}^{\infty} n \nu_\theta(n) = 2s \theta
\ee
The discrete generalized KMP process is self-dual \cite{carinci2013duality} with self-duality functions given by
\be\label{diskmpdual}
D_F(\xi,\eta)= \prod_i \frac{\eta_i!}{(\eta_i-\xi_i)!} \frac{\Gamma(2s)}{\Gamma(2s + \xi_i)}
\ee
More precisely, we have
\be\label{selfdualdiskmp}
\E_\eta \big(D_F(\xi, \eta(t))\big)= \E_\xi\big( D_F(\xi(t), \eta)\big)
\ee
The subscript ``$F$'' is added to the duality function  $D_F$ to recall that,
for a given $\xi \in \N^V$, the expectation of the duality function
w.r.t. a measure on the $\eta$ variables gives essentially {\em the multivariate factorial
moments} (up to the factors  $\frac{\Gamma(2s)}{\Gamma(2s + \xi_i)}$).
In particular the relation between the self-duality polynomials  and product measures with marginals \eqref{discgam}
reads
\be\label{intdual}
\int D_F (\xi,\eta) \bigotimes_{i\in V} \nu_{\theta_i}[d\eta]= \prod_{i\in V} \theta_i^{\xi_i}
\ee
This equality completely characterizes the product measure $\bigotimes_{i\in V} \nu_{\theta_i}$
via its factorial moments.

The hidden parameter model associated to the discrete generalized KMP is a process on $[0,\infty)^V$ which will determine the evolution of the parameters $\theta=(\theta_i, i\in V$) of product measures of the type
$\bigotimes_i \nu_{\theta_i}$.
The process is defined in the spirit of \cite{deMasi} via its generator
\be\label{diskmphid}
\loc f(\theta)= \sum_{i,j} p(i,j) \loc_{ij} f(\theta)
\ee
where the single edge generator $\loc_{ij}$ acts on the variables $\theta_i,\theta_j$ as follows
\be\label{}
\loc_{ij} f(x,y)= \E \Big(f(xB+y(1-B), xB+y(1-B)) - f(x,y)\Big)
\ee
where $B$ has a Beta distribution with parameters $(2s, 2s)$, i.e., it has density \eqref{beta}.
More explicitly we have
\be\label{diskmphidsedge}
\loc_{ij} f(x,y)= \int_0^1\big(f( xu + y(1-u), xu + y (1-u))- f(x,y)\big) \Beta (2s,2s) [du]
\ee

We then have the following duality result.
\bp\label{opi}
The  discrete generalized KMP process with generator \eqref{diskmp} is dual to the hidden parameter model with
generator \eqref{diskmphid} with duality function
\be
D(\xi, \theta) = \prod_{i\in V} \theta_i^{\xi_i}
\ee
\ep
\bpr
We act with the generator $L_{ij}$ in \eqref{diskmpsedge} on the $\xi$ variables
and obtain, using the binomial formula
\beq\label{byio}
L_{ij} \theta_i^{\xi_i}\theta_j^{\xi_j}
&=& \E\left(\theta_i^{X}\theta_j^{\xi_j+\xi_i-X}\right)- \theta_i^{\xi_i}\theta_j^{\xi_j}
\nonumber\\
&=&
\sum_{k=0}^{\xi_i+\xi_j} {\xi_i+\xi_j\choose k}\int_0^1 p^k (1-p)^{ \xi_i+\xi_j-k} \theta_i^k \theta_j^{\xi_i+\xi_j-k} \Beta(2s,2s)[dp]  - \theta_i^{\xi_i}\theta_j^{\xi_j}
\nonumber\\
&=&
\int_0^1 (p\theta_i +(1-p)\theta_j)^{\xi_i+\xi_j} \Beta(2s,2s)[dp] - \theta_i^{\xi_i}\theta_j^{\xi_j}
\eeq
This is now clearly the same as acting with the generator $\loc_{ij}$ in \eqref{diskmphidsedge} on the $\theta$ variables.
\epr

We can then state a result on the evolution of product measures of the type
$\bigotimes_{i\in V} \nu_{\theta_i}$ under the  discrete generalized KMP model.
\bt\label{prodth}
Consider the  discrete generalized KMP model  with generator \eqref{diskmp} and start it from a product measure
$
\bigotimes_{i\in V} \nu_{\theta_i}.
$
Denote by
$
\big(\bigotimes_{i\in V} \nu_{\theta_i}\big) S(t)
$
the evolved measure at time $t>0$, where $(S(t))_{t\ge 0}$ is the semigroup.  Then we have
\be
\big(\bigotimes_{i\in V }\nu_{\theta_i}\big) S(t)= \E_{\theta} \big(\bigotimes_{i\in V} \nu_{\theta_i(t)}\big)
\ee
where $\E_\theta$ denotes the expectation in the hidden parameter model with
generator \eqref{diskmphid} initialized
from the configuration $\theta$.
As a consequence, the set of mixtures
\[
\int \big( \bigotimes_{i\in V} \nu_{\theta_i}\big) \Xi [d\theta]
\]
is closed under the evolution of the  discrete generalized KMP model.
\et
\bpr
The proof uses self-duality of the discrete generalized KMP process (stated in \eqref{selfdualdiskmp}), and the duality between discrete generalized KMP and
the hidden parameter model (Proposition \ref{opi}). As a consequence of the identity \eqref{intdual} we  obtain the following series of equality:
\beq
\int D_F(\xi,\eta) \big(\bigotimes_{i\in V} \nu_{\theta_i}\big) S(t)[d\eta]
&=&
\int \E_\eta \big(D_F(\xi,\eta(t))\big)\big( \bigotimes_{i\in V} \nu_{\theta_i}\big)[d\eta]
\nonumber\\
&=& \int \E_\xi \big(D_F(\xi(t),\eta)\big) \big(\bigotimes_{i\in V} \nu_{\theta_i}\big)[d\eta]
\nonumber\\
&=&
\E_\xi \Big(\prod_{i\in V} \theta_i^{\xi_i(t)}\Big)
\nonumber\\
&=&
\E_\theta \Big(\prod_{i\in V} \theta_i(t)^{\xi_i}\Big)
\nonumber\\
&=& \E_\theta \int D_F (\xi,\eta) \big(\bigotimes_{i\in V} \nu_{\theta_i(t)}\big)[d\eta]
\eeq
Here in the second equality we used self-duality of the discrete generalized KMP process and in the fourth equality we used Proposition \ref{opi}.
The proof is then completed by observing that the functions $\eta\to D(\xi,\eta)$ are measure determining.
\epr

The result of Theorem \ref{prodth} can be reformulated as an intertwining result between the hidden parameter process and the discrete generalized KMP process.
We say that two Markov  processes with semigroups $(S(t), t\geq 0)$ and $(\caS(t), t\geq 0)$ are intertwined with intertwiner $\caG$ if for all $t\geq 0$
\be\label{intertwi}
\caG S(t) = \caS(t) \caG 
\ee
In Theorem \ref{prodth} we have obtained
\be\label{intopi}
\int S(t) f(\eta) \big(\bigotimes_{i\in V} \nu_{\theta_i}\big) [d\eta]= \caS(t) \int f(\eta) \big(\bigotimes_{i\in V} \nu_{\theta_i}\big) [d\eta]
\ee
where $S(t)$ is the semigroup of the discrete generalized KMP process and where $\caS(t)$ is the semigroup of
the hidden parameter model.
Therefore, if we define for a function $f:\N^V\to\R$ the ``discrete-gamma'' intertwiner
\[
\caG f(\theta)= \int f(\eta)\big( \bigotimes_{i\in V} \nu_{\theta_i}\big) [d\eta]
\]
where we implicitly assumed that $f$ is integrable w.r.t. $\bigotimes_{i\in V} \nu_{\theta_i}$,
then \eqref{intopi} reads
\[
\caG (S(t) f)= \caS(t) (\caG f)
\]
which is exactly intertwining between the hidden parameter process and the discrete KMP process.
\subsection{Continuous generalized KMP}
The continuous generalized KMP process with parameter $2s>0$ is a process on $[0,\infty)^V$ and is defined via the generator
\be\label{ckmp}
L f(\zeta)
= \sum_{i,j\in V} p(i,j) L_{ij} f(\zeta)
\ee
where the single edge generator $L_{ij}$ works on the variables $\eta_i,\eta_j$ as follows
\be\label{ckmp-edge}
L_{ij} f(x,y)= \E \big( f(X(x+y), (1-X)(x+y)) - f(x,y) \big)\,.
\ee
Here $X$ is $\Beta(2s,2s)$ distributed random variable.

The reversible measures of the continuous generalized KMP process are products of Gamma distribution
with parameters $(\theta,2s)$, 
where $\theta$ is the scale parameter and where $2s$ is the shape parameter, i.e.  the marginals are given by
\be\label{gam}
\nu_\theta[dx]=\frac{x^{2s-1}}{\theta^{2s}\Gamma(2s)}e^{-x/\theta} \ dx
\ee
The continuous and discrete generalized KMP processes are dual \cite{kmp,carinci2013duality} with duality function
\[
D_m(\xi, \zeta)= \prod_{i\in V} \zeta_i^{\xi_i} \frac{\Gamma(2s)}{\Gamma(2s+\xi_i)}
\]
The subscript ``$m$'' is added to the duality function  $D_m$ to recall that,
for a given $\xi \in \N^V$, the expectation of the duality function
w.r.t. a measure on the $\zeta$ variables gives essentially {\em the multivariate
moments} (up to the factors  $\frac{\Gamma(2s)}{\Gamma(2s + \xi_i)}$).
In particular the relation between the duality functions and product measures with marginals \eqref{gam} reads
\be\label{contchar}
\int  D_m(\xi, \zeta) \big(\bigotimes_{i\in V} \nu_{\theta_i} \big)[d\zeta]= \prod_{i\in V} \theta_i^{\xi_i}
\ee
This equality completely characterizes the product measure $\bigotimes_{i\in V} \nu_{\theta_i}$
via its  moments.

The main result on the evolution of product measures of the type
$\bigotimes_{i\in V} \nu_{\theta_i}$ under the continuous generalized KMP model
is stated in the following theorem.
\bt\label{prodexth}
Start the continuous generalized KMP model with generator \eqref{ckmp} from a product measure
$\bigotimes_{i\in V}\nu_{\theta_i} [d\zeta]$.
Then at time $t>0$ we have the measure
\be
\big(\bigotimes_{i\in V}\nu_{\theta_i}\big) S(t) [d\zeta]= \E_\theta\big(\bigotimes_{i\in V} \nu_{\theta_i(t)} [d\zeta]\big)
\ee
where $\{\theta(t), t\geq 0\}$ is the hidden parameter model with
generator \eqref{diskmphid} initialized from the configuration $\theta$.
\et
\bpr
We use the duality between the continuous and discrete generalized KMP model, combined the duality between the discrete KMP model and the hidden parameter model.
We then obtain
\beq
\int D_m(\xi,\zeta)\big(\bigotimes_{i\in V} \nu_{\theta_i}\big) S(t)[d\zeta]
&=&
\int \E_\eta \big(D_m (\xi,\zeta(t))\big) \big(\bigotimes_{i\in V} \nu_{\theta_i}\big)[d\zeta]
\nonumber\\
&=& \int \E_\xi \big(D_m (\xi(t),\zeta)\big)\big( \bigotimes_{i\in V} \nu_{\theta_i}\big)[d\zeta]
\nonumber\\
&=&
\E_\xi \Big(\prod_i \theta_i^{\xi_i(t)}\Big)
\nonumber\\
&=&
\E_\theta \Big(\prod_i \theta_i(t)^{\xi_i}\Big)
\nonumber\\
&=& \E_\theta \int D_m (\xi,\zeta)  \big(\bigotimes_{i\in V} \nu_{\theta_i(t)}\big)[d\zeta]
\eeq
We then conclude by observing that the functions $\zeta\to D_m (\xi,\zeta)$ are measure determining.
\epr

We then have the analogous result of Proposition \ref{opi} in the setting of the continuous generalized KMP process.
\bp\label{copi}
The continuous generalized KMP process with generator \eqref{ckmp} and the hidden parameter model  with
generator \eqref{diskmphid} are dual with duality function
\be\label{lap}
D(\theta, \zeta)= \prod_{i\in V} e^{\theta_i \zeta_i}
\ee
\ep
\bpr
It suffices to prove the duality for the single edge generators. Acting with the single edge generator of the continuous generalizd KMP model on the $\zeta$ variables gives
\beq
&&L_{ij}e^{\theta_i \zeta_i}e^{\theta_j \zeta_j}
\nonumber\\
&=& \int_0^1 \left(e^{\theta_i u(\zeta_i+\zeta_j) +\theta_j (1-u) (\zeta_i+\zeta_j)}- e^{\theta_i \zeta_i}e^{\theta_j \zeta_j}\right) \Beta(2s,2s)[du]
\nonumber\\
&=&
\int_0^1 \left(e^{(u\theta_i + (1-u) \theta_j)\zeta_i + (u\theta_i + (1-u) \theta_j)\zeta_j} -e^{\theta_i \zeta_i}e^{\theta_j \zeta_j}\right) \Beta(2s,2s)[du]
\eeq
which is recognized as the action of the generator $\loc_{ij}$ in \eqref{diskmphidsedge} on the $\theta$ variables.
\epr
\br\label{rempoiss}
Notice that we can find the duality function between continuous generalized KMP and the hidden parameter model also via the generating function of the duality function between discrete generalized KMP and the hidden parameter model, i.e.,
\[
\sum_{n=0}^\infty\frac{\theta^n z^n}{n!}= e^{\theta z}
\]
Indeed,  the continuous and discrete generalized KMP model are intertwined via the intertwiner
\[
\Lambda f(z)=\sum_{n=0}^\infty f(n) \frac{z^n}{n!}
\]
More precisely denoting here by $L_{d}$ the generator of the discrete generalized KMP \eqref{diskmp}
and by  $L_{c}$ the generator of the continuous generalized KMP  \eqref{ckmp}, we have for $f:\N^V\to\R$
\[
\Lambda(L_{d}f)= L_{c}(\Lambda f)
\]
where with a small abuse of notation we denoted by $\Lambda$ the tensorization of $\Lambda$, i.e.,
the $\Lambda$ acting on all the variables $\eta_i$
\[
\Lambda f(\zeta)= \sum_{\eta\in \N^V} f(\eta) \frac{\zeta^\eta}{\eta!}
\]
where
\[
\frac{\zeta^\eta}{\eta!}= \prod_{i\in V}\frac{\zeta_i^{\eta_i}}{\eta_i!}
\]
\er
Also here, we can reformulate Theorem \ref{prodexth} as an intertwining result. Indeed, by considering the Gamma distribution in \eqref{gam} and by defining
the ``Gamma'' intertwiner
\[
\caG f(\theta) = \int f(\zeta) \bigotimes_{i \in V}\nu_{\theta_i} [d\zeta]
\]
it follows that Theorem \ref{prodexth} can be read as an intertwining between the hidden parameter process and the continuous generalized KMP process, with intertwiner $\caG$.

\subsection{Adding  driving}
We will discuss the adding of driving for the continuous generalized KMP model only.
The results for the generalized discrete KMP model are completely analogous.

We start by describing the generator modelling the coupling to a reservoir.
It is a  generator that acts on a single variable $x\in\R$ as follows
\be\label{boundi}
L_{\theta^*} f(x)= \E(f((x+Y)B)- f(x))
\ee
where $\E$ denotes expectation over the two independent random variables $U,Y$
and where $Y$ is distributed as $\nu_{\theta^*}$ (Gamma distribution) and $B$ is $\Beta(2s,2s)$ distributed.
Thus the action of the boundary site reservoir generator is similar to the bulk edge generator,
in the sense that the redistribution of energies between the site and the reservoir
occurs via a Beta random variable; however now the energy of the ``extra site''
representing the reservoir is sampled from a Gamma distribution with
mean $2s \theta^*$, which is exactly the marginal of the invariant distribution of the model without reservoirs.
Reservoirs of this form were introduced originally in the setting of the KMP model (corresponding to $2s=1$) in \cite{bert} and are different from the reservoirs in the original model \cite{kmp}.

The corresponding boundary generator of the hidden parameter model is
\be\label{boundihidd}
\loc_{\theta^*} f(\theta)= \int_0^1 \Big(f((1-u) \theta+ u \theta^*) - f(\theta)\Big) \Beta(2s,2s)[du]
\ee
which can be viewed as having an ``extra site'' from which always
 the value $\theta_*$ is imported.

We then have the following intertwining result.
\bl\label{interlem}
For a function $f:[0,+\infty)\to \mathbb{R}$ which is integrable with respect to
the Gamma distribution $\nu_{\theta}$ define the intertwiner
$$
\caG f(\theta)= \frac{1}{\Gamma(2s) \theta^{2s}}\int_0^{\infty} f(x) {x^{2s-1}}{e^{-x/\theta}} dx
$$
Then the boundary generator of the continuous generalized KMP process \eqref{boundi}
and the boundary generator of the hidden parameter model \eqref{boundihidd}
are intertwined as
\be\label{boundin}
\caG L_{\theta^*}  = \loc_{\theta^*} \caG
\ee
\el
\bpr
For simplicity we prove the case $2s=1$, the general case is obtained with a similar proof.
We have
\begin{eqnarray}
(\caG L_{\theta^*} f)(\theta)
& = &
\int_{0}^\infty dx \frac{e^{-x/\theta}}{\theta} L_{\theta^*} f(x) \nonumber\\
& = &
\int_0^\infty dx \frac{e^{-x/\theta}}{\theta} \int_0^{\infty} dy \frac{e^{-y/\theta^*}}{\theta^*} \int_0^1 du \Big(f((x+y)u) - f(x) \Big) \nonumber\\
\end{eqnarray}
and we also have
\begin{eqnarray}
(\loc_{\theta^*} \caG f)(\theta)
& = &
\int_{0}^1 du \Big(  \caG f((1-u) \theta+ u \theta^*) -  \caG f(\theta)\Big)  \nonumber\\
& = &
\int_{0}^1 du \Big(  \int_{0}^{\infty} dx \frac{e^{-\frac{x}{(1-u) \theta+ u \theta^*}}}{(1-u) \theta+ u \theta^*}  f(x) -  \int_{0}^{\infty} dx \frac{e^{-\frac{x}{\theta}}}{\theta}  f(x)\Big)
\end{eqnarray}
It suffices to see \eqref{boundin} for the functions $f_n(x)= x^n/n!$ (for all $n\in\N$).
From the previous two equations this in turn reduces to proving the following identity
\beq\label{balanka}
\int_0^\infty dx \frac{e^{-x/\theta}}{\theta} \int_0^{\infty} dy \frac{e^{-y/\theta^*}}{\theta^*} \int_0^1 du \frac{((x+y)u)^n}{n!}
=
\int_{0}^1 du  \int_{0}^{\infty} dx \frac{e^{-\frac{x}{(1-u) \theta+ u \theta^*}}}{(1-u) \theta+ u \theta^*}  \frac{x^n}{n!}
\eeq
The rhs of \eqref{balanka} equals
\be\label{ri}
\int_{0}^1 du  \int_{0}^{\infty} dx \frac{e^{-\frac{x}{(1-u) \theta+ u \theta^*}}}{(1-u) \theta+ u \theta^*}  \frac{x^n}{n!}
= \int_0^1 du (u\theta_*+ (1-u)\theta)^n= \frac1{n+1}\sum_{k=0}^n (\theta_*)^k \theta^{n-k}
\ee
where we used the identity
\[
\int_0^1 u^k (1-u)^{n-k} du= \frac{k!(n-k)!}{(n+1)!}
\]
combined with
$
\int_0^\infty \frac{x^n}{n!} \frac{e^{-x/\theta}}{\theta} dx = \theta^n
$.
The lhs of \eqref{balanka} equals
\beq
&&\int_0^\infty dx \frac{e^{-x/\theta}}{\theta} \int_0^{\infty} dy \frac{e^{-y/\theta^*}}{\theta^*} \int_0^1 du \frac{((x+y)u)^n}{n!}
\nonumber\\
&=&
\int_0^\infty dx \frac{e^{-x/\theta}}{\theta} \int_0^{\infty} dy \frac{e^{-y/\theta^*}}{\theta^*} \frac{1}{(n+1)} \sum_{k=0}^n   \frac{x^k}{k!} \frac{y^{n-k}}{(n-k)!}
\nonumber\\
&=&
\frac{1}{n+1}\sum_{k=0}^n \theta^k (\theta_*)^{n-k}
\eeq
\epr

To define the general boundary driven model, we associate reservoirs with parameters $\theta^*_i$ at site $i\in V$
and the generator of the boundary driven  continuous generalized KMP process is then given by
\be\label{bdriv}
L f(\zeta)
= \sum_{i,j\in V} p(i,j) L_{ij} f(\zeta) + \sum_{ i\in V} c(i)  L_{\theta^*_i}f(\zeta)
\ee
where $L_{i,j}$ is read in \eqref{ckmp-edge} and  $L_{\theta^*_i}$ is defined in \eqref{boundi}.
As a consequence of Theorem \ref{interlem} and of the intertwining result of Lemma \ref{prodth},
we then have the following propagation of mixtures of products of Gamma distributions.
\bt\label{hiddenKMP}
Consider the  driven continuous generalized KMP model with generator \eqref{bdriv}. 
Then we have the following.
\begin{itemize}
\item[a)]  If we  start the process from a product measure of the form
$
\bigotimes_{i\in V} \nu_{\theta_i},
$
then at time $t>0$ the distribution is given by
\[
\E_{\theta} \big(\bigotimes_{i\in V} \nu_{\theta_i(t)}\big)
\]
where the process $(\theta_i(t), i\in V, t\geq 0)$ evolves according to the generator
\be\label{opibound}
\sum_{i,j\in V}p(i,j)\loc_{ij} + \sum_{i\in V} c(i) \loc_{\theta^*_i}
\ee
\item[b)] The  driven generalized KMP process converges to a unique stationary measure which reads
\[
\int \bigotimes_{i\in V} \nu_{\theta_i} \Xi[d\theta]
\]
where the mixture measure $\Xi$ is the unique stationary measures of the associated hidden parameter model, with generator \eqref{opibound}.
\item[c)]
In particular if all the reservoir parameters are equal to a fixed value, i.e. $\theta^*_i =\theta^*$ for all $i\in V$,
then this unique stationary measure is given by $\bigotimes_{i\in V} \nu_{\theta^*}$ and is  also reversible.
\end{itemize}
\et
\section{Generalized harmonic models}
\label{section3}

In this section we consider the generalized discrete harmonic model \cite{frass1}  and the associated generalized continuum harmonic model
(also called integrable heat conduction model in \cite{frass2}).
The aim here is to prove the existence of a hidden parameter model and to derive conclusions
from it about the nature of the stationary measures in the one-dimensional boundary driven set-up.
Contrary to the KMP model, the invariant measure of the hidden parameter model on the chain with left and right boundary reservoirs can be obtained explicitly. The main reason is a hidden Markovian structure of the hidden parameter model, see section \ref{invhid} and \ref{invhidgen} below for details. This hidden Markovian structure can be seen as the probabilistic counterpart of the integrability of this model, which was used in previous works \cite{frass1}, \cite{frass2}, \cite{red} to obtain the non-equilibrium steady state on the chain.

\subsection{Mass redistribution models}
In order to introduce  the harmonic models, let us first consider the following general class of generators (see also \cite{capanna2024class}) acting on two variables $y_1, y_2\geq 0$, and parametrized by a positive measure $M$ on the interval $[0,1]$.
\beq\label{gengen}
L_{12}f(y_1,y_2)
&=& L^{M}_{12}f(y_1, y_2)
\\
&=& \int_0^1 M(du) \left[(f(y_1-uy_1, y_2+ u y_1)
+f(y_1+uy_2, y_2- u y_2) -2f(y_1, y_2))\right]\nonumber
\eeq
In this process, with rate $M(du)$, a fraction of mass is taken away from one of the two sites and given to the other site.
Notice that in these models, different from the KMP model, only a fraction of the mass of {\em one site} is moved to the other site (rather than a fraction of the total mass of the two sites).

In order to introduce the associated hidden parameter model, we consider the following generator acting on two variables $\theta_1,\theta_2\geq 0$
\be\label{genopi}
\loc_{12} f(\theta_1,\theta_2)= \int_0^1 M(du) [ f(\theta_1 (1-u)+ u \theta_2, \theta_2)+ f(\theta_1, u\theta_1+ (1-u)\theta_2)- 2f(\theta_1,\theta_2)]
\ee
We see that, contrary to the hidden parameter model for the generalized KMP processes, here the parameters (or ``local temperatures'') $\theta_1,\theta_2$ are replaced by convex combinations only at one of the two sites, leaving the parameter at the other site untouched. We have the following duality result.
\bp\label{mdual}
The process with generator $L_{12}$ in  \eqref{gengen} is dual to the process with generator $\loc_{12}$ in \eqref{genopi}
with duality function
\[
D_c(\theta_1,\theta_2;y_1,y_2)= e^{\theta_1y_1+\theta_2y_2}
\]
\ep
\bpr
This follows from the simple observation
\[
e^{\theta_1 (y_1-uy_1) + \theta_2 (y_2+u y_1)}= e^{(\theta_1(1-u) +\theta_2 u) y_1 +\theta_2 y_2}
\]
and the similar equality obtained by interchanging  the sub-indices $1$ and $2$.
\epr

To understand associated intertwined discrete models, let us consider the Poisson intertwiner
between functions $f: \N^2\to\R$ and functions $f: [0,\infty)^2\to\R$
\be\label{poiss}
\Lambda_{12} f(y_1, y_2)=\sum_{k_1, k_2 \in \N} f(k_1, k_2)\frac{y_1^{k_1}}{k_1!} \frac{y_2^{k_2}}{k_2!}
\ee
Now we consider discrete models of mass redistribution, i.e., Markov processes on $\N^2$ depending on a positive measure
$\caM(k,n), k\in\N, n\in\N$ with support $\{(k,n): k\leq n\}$. 
The discrete models are then defined via their generator acting on functions $f:\N^2\to\R$ as follows
\beq\label{disgengen}
Lf(n_1, n_2) &=& L^{\caM}_{12}f(n_1,n_2)
\nonumber\\
&=&\sum_{k=1}^{n_1} \caM(k,n_1)(f(n_1-k, n_2+k)- f(n_1,n_2))
\nonumber\\
&+& \sum_{k=1}^{n_2} \caM(k,n_2)(f(n_1+k, n_2-k)- f(n_1,n_2))
\eeq
We say that the discrete model \eqref{disgengen} {\em is associated to the the continuum model} \eqref{gengen} if it is Poisson intertwined with it, i.e., if
\be\label{poisin}
\Lambda_{12} (L^{\caM}_{12})= L^M_{12}(\Lambda_{12})
\ee
Then we have the following lemma relating dualities of the continuous models to dualities of the associated discrete models.
\bl\label{poislem}
If \eqref{poisin} holds, then the process with generator $L^{\caM}_{12}$ is dual to the process with generator \eqref{genopi} with duality function
\[
D_d(\theta_1,\theta_2; n_1, n_2)= \theta_1^{n_1}\theta_2^{n_2}
\]
\el
\bpr
This follows by the following two facts:  i) the Poissonian generating function applied to $\theta_1^{n_1}\theta_2^{n_2}$ equals
$e^{\theta_1 y_1+\theta_2y_2}$ (cf. Remark \ref{rempoiss}); ii) the duality between the discrete process with generator $L_{12}^\caM$
with duality functions $D_d(\theta_1,\theta_2; n_1,n_2)$ and the process with generator \eqref{genopi} is equivalent with duality between the continuous process
with generator $L_{12}^M$ and the process with generator \eqref{genopi} with duality
function
\[
D_c(\theta_1,\theta_2; n_1,n_2) = \sum_{n_1, n_2}D_d(\theta_1,\theta_2; n_1,n_2) \frac{y_1^{n_1} y_2^{n_2}}{n_1! n_2!}
\]
See e.g. \cite{book} for a proof of this equivalence.
The duality for continuum models of Proposition \eqref{mdual} therefore implies automatically the duality for discrete models which are Poisson intertwined.
\epr
\subsection{The harmonic models}
For the simplest version of the continuous harmonic model, we have
$M(du)=\frac{1}{u} du$ and for the associated discrete model $\caM(k,n)= \frac{1}{k} I(1\leq k\leq n)$.
We first consider the model on a general graph with vertex set $V$ and with edge weights $p(i,j)$
and define the generator on functions $f:[0,\infty)^V\to\R$ as follows
\be\label{charmhalf}
L f(\zeta)= \sum_{i,j\in V} p(i,j)
\int_0^1 \frac{du}{u} \left((f(\zeta- u\zeta_i\delta_i + u\zeta_i \delta_j)-f(\zeta))
+ (f(\zeta- u\zeta_j\delta_j + u\zeta_j \delta_i)-f(\zeta))\right)
\ee
Here $\delta_i$ denotes the configuration with unit mass at site $i$ and zero mass everywhere else.
The process corresponding to the generator \eqref{charmhalf} will be called the continuous harmonic process.
Its reversible product measures are products of exponentials with identical scale parameters, i.e., 
with marginals
\be\label{expo}
\nu_\theta (dx)= \frac1{\theta}e^{-x/\theta}
\ee
with expectation $\theta >0$.
The associated discrete model is then defined via its generator acting in functions
$f:\N^V\to\R$:
\be\label{harmhalf}
L f(\eta)= \sum_{i,j\in V} p(i,j)
\left(\sum_{k=1}^{\eta_i}\frac1k (f(\eta- k\delta_i + k\eta_j \delta_j)-f(\eta))
+\sum_{k=1}^{\eta_j}\frac1k(f(\eta- k\delta_j + k\delta_i)-f(\zeta))\right)
\ee
We call the corresponding process the discrete harmonic process.
Its reversible product measures are products of geometric random variables
with marginals
\be\label{geom}
\nu_\theta (n)= \left(\frac{\theta}{1+\theta}\right)^n \left(\frac{1}{1+\theta}\right)
\ee
with mean $\theta$.
Finally, the corresponding hidden parameter model is defined via its generator
acting on $f:[0,\infty)^V\to\R$:
\be\label{opiharm}
L f(\theta)= \sum_{i,j\in V} p(i,j)
\int_0^1 \frac{du}{u} \left((f(\theta- u\theta_i\delta_i + u\theta_j \delta_i)-f(\theta))
+ (f(\theta- u\theta_j\delta_j + u\theta_i \delta_j)-f(\theta))\right)
\ee
This generator was also considered in the literature of integrable systems, see for instance Eq. (2.3.3) in  \cite{derkachov1999baxter} where it appears as a representation of the integrable XXX spin chain, and section 2.3 in \cite{frassek2020non} where a connection between the generator \eqref{opiharm} and the continuous harmonic
generator was pointed out. 
\bt\label{grandhalfdualitythm}
\begin{itemize}
We have the following duality and intertwining relations:
\item[a)] The discrete harmonic model is self-dual with self-duality function
\be
D_F(\xi, \eta)= \prod_{i\in V} {\eta_i \choose \xi_i}
\ee
\item[b)] The discrete and continuous harmonic models are dual with duality function
\be
D_m(\xi, \zeta)= \prod_{i\in V} \frac{\zeta_i^{\xi_i}}{\xi_i!}
\ee
\item[c)] The discrete and continuous harmonic model are Poisson intertwined.
\item[d)] The continuous harmonic process and the hidden parameter model are dual with duality function
\be
D_c(\theta,\zeta)= \prod_{i\in V} e^{\theta_i\zeta_i}
\ee
\item[e)] The discrete harmonic process and the hidden parameter model are dual with duality function
\be
D_d(\theta,\eta)= \prod_{i\in V} \theta_i^{\eta_i}
\ee
\end{itemize}
\et
\bpr
See \cite{frass1}, \cite{frass2} for the statements a) up to c). From c) it follows that d) and e) are equivalent via Lemma
\ref{poislem}, and d) follows from Proposition \ref{mdual}.
\epr

We can then turn the duality results into a result on propagation of mixtures of product measures,
or equivalently into an intertwining result.
\bt\label{mixthmharm}
The following results hold.
\begin{itemize}
\item[a)] Start the discrete harmonic process with generator \eqref{harmhalf} from a product measure with geometric marginals
$\bigotimes_{i\in V} \nu_{\theta_i}$, where $\nu_\theta$ is as in \eqref{geom}.
Then at time $t>0$ the distribution is equal to
\be\label{colanko}
\big(\bigotimes_{i\in V} \nu_{\theta_i}\big) S(t)= \E_\theta\big(\bigotimes_{i\in V} \nu_{\theta_i(t)}\big)
\ee
where $(\theta(t), t\geq 0)$ is the hidden parameter process with generator \eqref{opiharm}, and $\E_\theta$ denotes expectation in this process starting from $\theta$.
Equivalently, considering the ``geometric'' intertwiner of an integrable function $f:\N^V\to\R$
\[
\caG f(\theta)= \int f(\eta)\big( \bigotimes_{i\in V} \nu_{\theta_i}\big) [d\eta]
\]
we have
\[
\caG (L f)= \caL(\caG f)
\]
which is the intertwining between the generator $\caL$ of the hidden parameter process and the generator $L$ of the discrete harmonic process.

\item[b)] Start the continuous harmonic process  with generator \eqref{charmhalf}  from a product measure with exponential marginals
$\otimes \nu_{\theta_i}$, where $\nu_\theta$ is as in \eqref{expo}. Then at time $t>0$ the distribution is equal to
\be\label{biliflaki}
\big(\bigotimes_{i\in V} \nu_{\theta_i}\big) S(t)= \E_\theta\big(\bigotimes_{i\in V} \nu_{\theta_i(t)}\big)
\ee
where $\theta(t), t\geq 0$ is the hidden parameter process with generator \eqref{opiharm}, and $\E_\theta$ denotes expectation in this process starting from $\theta$.
Equivalently, considering the ``exponential'' intertwiner of an integrable function $f:[0,\infty)^V\to\R$
\[
\caG f(\theta)= \int f(\zeta)\big( \bigotimes_{i\in V} \nu_{\theta_i}\big) [d\zeta]
\]
we have
\[
\caG (L f)= \caL(\caG f)
\]
which is the intertwining between the generator $\caL$ of the hidden parameter process and the generator $L$ of the continuous harmonic process.
\end{itemize}
\et
\bpr
We will prove \eqref{biliflaki}.
The proof of \eqref{colanko} is analogous, replacing exponentials by geometric distributions.

The duality functions between the continuous and discrete harmonic model are given by
\[
D_m(\xi, \zeta)= \prod_i \frac{\zeta_i^{\xi_i}}{\xi_i!}
\]
Then we obtain,
\beq
\int  D_m(\xi, \zeta) \big(\bigotimes_{i\in V}\nu_{\theta_i}[d\zeta_i]\big) S(t) &=& \int \E_\zeta \big(D_m(\xi, \zeta(t))\big) \big(\bigotimes_{i\in V} \nu_{\theta_i}[d\zeta_i]\big)
\nonumber\\
&=&\int \E_\xi\big( D_m(\xi(t), \zeta)\big) \big(\bigotimes_{i\in V} \nu_{\theta_i}[d\zeta_i]\big)
\nonumber\\
&=&
\E_\xi \Big(\prod_{i\in V} \theta_i^{\xi_i(t)}\Big)
\nonumber\\
&=&
\E_\theta \Big(\prod_{i\in V} \theta_i(t)^{\xi_i}\Big)
\nonumber\\
&=&
\E_\theta\Big(\int  D_m(\xi, \zeta) \big(\bigotimes_{i\in V} \nu_{\theta_i(t)}[d\zeta_i]\big)\Big)
\eeq
Here we used duality between the continuous and the discrete model in the second equality, and
duality between the discrete model and the hidden parameter model in the third equality.
We can then conclude \eqref{biliflaki} because the functions $\zeta\to D(\xi, \zeta)$ are measure determining.
\epr

\subsection{Boundary reservoirs}
We now discuss the intertwining of
the boundary generator of the continuous harmonic model. This reads \cite{frass2}
\be
\label{genboundcharm}
L_{\theta^*}f(x) = \int_{0}^1 \frac{du}{u} (f(x(1-u)) - f(x)) + \int_{0}^{\infty} \frac{du}{u} e^{-u} (f(x+u \theta^*) - f(x))
\ee
This generator is reversible w.r.t. the exponential distribution with mean $\theta^* >0$.
The corresponding boundary hidden parameter generator has the same structure of the boundary hidden parameter KMP generator, the main difference being that the uniform measure of the KMP model is here replaced by the measure $du/u$. It reads
\be\label{genboundcharmtheta}
\loc_{\theta^*} f(\theta)= \int_0^1 \Big(f((1-u) \theta+ u \theta^*) - f(\theta)\Big) \frac{du}{u}
\ee

We then have the following intertwining result.
\bl\label{interlem22}
For a function $f:[0,+\infty)\to \mathbb{R}$ which is integrable with respect to
the Exponential distribution $\nu_{\theta}$ define the intertwiner
$$
\caG f(\theta)= \int_0^{\infty} f(x) \frac{e^{-x/\theta}}{\theta} dx
$$
The boundary generator of the continuous harmonic process \eqref{genboundcharm}
and the boundary generator of the hidden parameter model \eqref{genboundcharmtheta}
are intertwined as
\be\label{boundin22}
\caG L_{\theta^*}  = \loc_{\theta^*} \caG
\ee
\el
\bpr
We have
\begin{eqnarray}
\label{use1}
(\caG L_{\theta^*} f)(\theta)
& = &
\int_{0}^\infty dx \frac{e^{-x/\theta}}{\theta} L_{\theta^*} f(x) \\
& = &
\int_0^\infty dx \frac{e^{-x/\theta}}{\theta} \left( \int_{0}^1 \frac{du}{u} (f(x(1-u)) - f(x)) + \int_{0}^{\infty}\frac{du}{u} e^{-u}(f(x+u \theta^*) - f(x))\right)
 \nonumber
\end{eqnarray}
and we also have
\begin{eqnarray}
\label{use2}
(\loc_{\theta^*} \caG f)(\theta)
& = &
\int_{0}^1 \frac{du}{u} \Big(  \caG f((1-u) \theta+ u \theta^*) -  \caG f(\theta)\Big)  \\
& = &
\int_{0}^1 \frac{du}{u} \Big(  \int_{0}^{\infty} dx \frac{e^{-\frac{x}{(1-u) \theta+ u \theta^*}}}{(1-u) \theta+ u \theta^*}  f(x) -  \int_{0}^{\infty} dx \frac{e^{-\frac{x}{\theta}}}{\theta}  f(x)\Big)
\nonumber
\end{eqnarray}
It suffices to see \eqref{boundin22} for the functions $f(x)= x^n/n!$ (for all $n\in\N$). I.e., we have to prove
\be\label{boundin222}
\caG L_{\theta^*} f = \loc_{\theta^*} \caG f
\ee
for those $f$.
Plugging this $f$ into \eqref{use1} we get
\begin{eqnarray}
(\caG L_{\theta^*} f)(\theta)
& = &
\int_0^\infty dx \frac{e^{-x/\theta}}{\theta} \left( \int_{0}^1 \frac{du}{u} \left(\frac{x^n(1-u)^n}{n!} - \frac{x^n}{n!}\right) +  \int_{0}^{\infty} \frac{du}{u} e^{-u} \left(\frac{(x+u \theta^*)^n}{n!} - \frac{x^n}{n!}\right)\right)
 \nonumber\\
 & = &
 \theta^n \int_{0}^1 \frac{du}{u} ((1-u)^n-1) + \sum_{k=1}^n \theta^{n-k} (\theta^*)^k \frac{1}{k!}\int_{0}^{\infty}\frac{du}{u} e^{-u}u^{k}
  \nonumber\\
 & = &
 \theta^n \int_{0}^1 \frac{du}{u} ((1-u)^n-1) + \sum_{k=1}^n \theta^{n-k} (\theta^*)^k \frac{1}{k}
 \end{eqnarray}
 Plugging $f(x) = x^n/n!$ into \eqref{use2} we get
 \begin{eqnarray}
(\loc_{\theta^*} \caG f)(\theta)
& = &
\int_{0}^1 \frac{du}{u} \Big(  \int_{0}^{\infty} dx \frac{e^{-\frac{x}{(1-u) \theta+ u \theta^*}}}{(1-u) \theta+ u \theta^*}  \frac{x^n}{n!} -  \int_{0}^{\infty} dx \frac{e^{-\frac{x}{\theta}}}{\theta}  \frac{x^n}{n!}\Big)
\nonumber\\
& = &
\int_{0}^1 \frac{du}{u} \Big(  ((1-u) \theta+ u \theta^*)^n -  \theta^n \Big)
\nonumber\\
 & = &
 \theta^n \int_{0}^1 \frac{du}{u} ((1-u)^n-1) + \sum_{k=1}^n \theta^{n-k} (\theta^*)^k \frac{1}{k}
\end{eqnarray}
This completes the proof
\epr

We define the generator of the boundary driven continuous harmonic process
\begin{eqnarray}
\label{bdriv22}
L f(\zeta)
& = & \sum_{i,j\in V} p(i,j)
\int_0^1 \frac1u \left((f(\zeta- u\zeta_i\delta_i + u\zeta_i \delta_j)-f(\zeta))
+ (f(\zeta- u\zeta_j\delta_j + u\zeta_j \delta_i)-f(\zeta))\right)
\nonumber \\
&+& \sum_{ i\in V} c(i) \left( \int_{0}^1 \frac{du}{u} (f(\zeta- u\zeta_i\delta_i) - f(\zeta)) +  \int_{0}^1 \frac{du}{u} e^{-u}(f(\zeta+u \theta_i^*\delta_i) - f(\zeta))\right)
\nonumber\\
&&
\end{eqnarray}
Here we associate reservoirs with parameters $\theta^*_i$ to site $i\in V$.

As a consequence of Theorem \ref{mixthmharm} and of the intertwining result of Lemma \ref{interlem22},
we then have for the boundary driven continuous harmonic process the following propagation of mixtures of product of exponential distribution.
\bt\label{contharminv}
Consider the resevoir driven continuous harmonic model with generator \eqref{bdriv22}. Start the model from a product measure of the form
$
\bigotimes_{i\in V} \nu_{\theta_i},
$
where $\nu_{\theta_i}(d\zeta_i)= \frac{e^{-x/\theta_i}}{\theta_i} d\zeta_i$.
Then at time $t>0$ the distribution is given by
\[
\E_{\theta} \big(\bigotimes_{i\in V} \nu_{\theta_i(t)}\big)
\]
where the process $(\theta_i(t), i\in V, t\geq 0)$ evolves according to the generator
\begin{eqnarray}\label{opiharm22}
L f(\theta)
&=&
\sum_{i,j\in V} p(i,j)
\int_0^1 \frac{du}{u} \left((f(\theta- u\theta_i\delta_i + u\theta_j \delta_i)-f(\theta))
+ (f(\theta- u\theta_j\delta_j + u\theta_i \delta_j)-f(\theta))\right)
\nonumber\\
&+& \sum_{i\in V} c(i) \int \Big(f(\theta- u\theta_i\delta_i + u\theta^*_i \delta_i) - f(\theta)\Big) \frac{du}{u}
\end{eqnarray}
When $t\to\infty$, the reservoir driven continuous harmonic process converges to a unique stationary measure which reads
\[
\int\Big(\bigotimes_{i\in V} \nu_{\theta_i}\Big) \Xi(\prod_{i\in V}d\theta_i)
\]
where the mixture measure $\Xi$ is the unique stationary measures of the associated hidden parameter model, with generator \eqref{opibound}.
\et

\subsection{Invariant measure of the single site hidden parameter model}

When we consider the harmonic model with a single site in contact with two reservoirs
with $\theta_L=0$ and $\theta_R=1$, the generator of the associated hidden parameter model 
reads as follows (cf.\ \eqref{genboundcharmtheta})
\beq\label{harmboundhid}
\loc^{0,1} f(\theta)&=& \int_0^1 \frac1u (f(\theta(1-u))-f(\theta))du
\nonumber\\
&+& \int_0^1 \frac1u (f(u+\theta(1-u))-f(\theta))du
\eeq
We then prove the following.
\bp\label{unifsingle}
The unique stationary distribution of the process with generator \eqref{harmboundhid} is the uniform distribution on $[0,1]$.
\ep
\bpr
To infer the stationarity of the uniform measure for the generator
\eqref{harmboundhid} it is convenient to 
consider the harmonic model on a single edge which is given by
\begin{eqnarray*}
L_{12} f(\zeta_1,\zeta_2)&=& \int_0^1\frac{du}{u} \left(f(\zeta_1(1-u), \zeta_2+\zeta_1 u)- f(\zeta_1,\zeta_2)\right)
\\
&+& \int_0^1\frac{du}{u} \left(f(\zeta_1 + \zeta_2 u, \zeta_2(1-u))- f(\zeta_1,\zeta_2)\right)
\end{eqnarray*}
In this model $\zeta_1+\zeta_2$ is conserved. Therefore, if we fix $\zeta_1+\zeta_2=1$ then, 
substituting $\zeta_2=1-\zeta_1$, we see that the action of the  generator $L_{12}$ on the $\zeta_1$ variable is exactly the same as the action of the generator \eqref{harmboundhid} on the $\theta$ variable.
We know that the reversible measures for the generator $L_{12}$ are
product measures with marginals exponentials with identical scale parameter, i.e., with joint density given by
\[
\frac1{\theta^2}e^{-\zeta_1/\theta} e^{-\zeta_2/\theta}
\]

As a consequence, considering two independent exponential random variables, the distribution of the first conditional to their sum being $1$ is invariant for
the generator \eqref{harmboundhid}. This is exactly the uniform distribution.
\epr

For general reservoir case $0\leq \theta_L\leq \theta_R$ the generator of the single-site hidden parameter model reads
\beq\label{harmboundhid22}
\loc^{\theta_L,\theta_R} f(\theta)&=& \int_0^1 \frac1u (f(\theta_L u +\theta(1-u))-f(\theta)) du
\nonumber\\
&+& \int_0^1 \frac1u (f(\theta_R u+\theta(1-u))-f(\theta)) du
\eeq
Using the change of variable $x \mapsto \theta_L+x(\theta_R-\theta_L)$,
Proposition \ref{unifsingle} implies that 
the uniform distribution on $[\theta_L,\theta_R]$ is
invariant for the process with generator $\loc^{\theta_L,\theta_R}$.

\subsection{The invariant measure of the hidden parameter model on the chain}\label{invhid}
In this section we consider the geometry of the chain $\{1, \ldots, N\}$ with boundary reservoirs
at left and right end. The hidden parameter model is then a model on the state space
$\Omega_N=[0,\infty)^{\{1,\ldots, N\}}$. It is parametrized by the left and right reservoir parameters $\theta_L, \theta_R$.
The generator of the hidden parameter model is  given by
\be\label{thetharm}
\loc f(\theta)= \loc_1^{\theta_L} f(\theta)+ \loc_N^{\theta_R} f(\theta) + \sum_{i=1}^N \loc_{i,i+1} f(\theta)
\ee
with boundary  single site generators 
\beq\label{folop1}
 \loc_1^{\theta_L} f(\theta) &=& \int_0^1 \frac{du}{u} \left(f(\theta-u\theta_1\delta_1+u\theta_L\delta_1)- f(\theta)\right)
\eeq
\beq\label{folopN}
 \loc_N^{\theta_R} f(\theta) &=& \int_0^1 \frac{du}{u} \left(f(\theta-u\theta_N\delta_N+u\theta_R\delta_N)- f(\theta)\right)
\eeq
and  with single edge generators 
\beq\label{folopbulk}
(\loc_{i,i+1} f)(\theta)= \int_{0}^1 \frac{du}{u} (f(\theta-u\theta_i\delta_i +u\theta_{i+1}\delta_{i})
+ f(\theta-u\theta_{i+1}\delta_{i+1} +u\theta_{i}\delta_{i+1}) -2f(\theta))
\eeq
In this subsection we identify the stationary measure of the boundary driven hidden parameter model on the chain.
As a consequence of Theorem \ref{contharminv} this yields also a full characterization of the
non-equilibrium steady state of the boundary driven continuous harmonic model on the chain (cf. also \cite{frass1,gab,red}).
\bt\label{bofi}
The invariant measure of the hidden parameter model with generator \eqref{thetharm} is the joint distribution of $(U_{1:1}, \ldots U_{N:N})$, the  order statistics of $N$ independent uniforms
on $[\theta_L,\theta_R]$.
As a consequence, the invariant measure of the boundary driven continuous harmonic model on the chain is
a mixture of product of exponential distributions with means $(U_{1:1}, \ldots U_{N:N})$.
\et
\bpr
We will prove the case $N=2, \theta_L=0, \theta_R=1$. As will turn out from the proof, the case $N=2$ is without loss of generality, via the Markovian structure of the joint distribution of order statistics, whereas the restriction $\theta_L=0, \theta_R=1$ can be generalized via elementary translation and scaling.
Let us call
$\Lambda^{2}_{0,1}$ the joint distribution of the order statistics of two independent uniforms on $[0,1]$. Let us call $\Lambda^1_{a,b}$ the distribution of one uniform on $[a,b]$.
Then we have the following conditional distributions
\[
\Lambda^{2}_{0,1}(d\theta_2|\theta_1=a)= \Lambda^1_{a,1} (d\theta_2), \  \Lambda^{2}_{0,1}(d\theta_1|\theta_2=b)= \Lambda^1_{0,b}(d\theta_1)
\]
So let us now consider the generator \eqref{thetharm} for $N=2, \theta_L=0, \theta_R=1$.
Then we want to prove that
\beq
\int \loc f(\theta_1,\theta_2) \Lambda^{2}_{0,1} (d\theta_1, d\theta_2)=0
\eeq
where
\beq\label{four}
\loc f(\theta_1,\theta_2) &=&
\int_0^1 \frac{du}{u}\left(f((1-u)\theta_1, \theta_2)-f(\theta_1, \theta_2)\right)
\nonumber\\
&+& \int_0^1 \frac{du}{u}\left(f((1-u)\theta_1+u\theta_2, \theta_2)-f(\theta_1, \theta_2)\right)
\nonumber\\
&+&
\int_0^1 \frac{du}{u}\left(f(\theta_1, (1-u)\theta_2+u\theta_1)-f(\theta_1, \theta_2)\right)
\nonumber\\
&+&
\int_0^1 \frac{du}{u}\left(f(\theta_1, u+(1-u)\theta_2)-f(\theta_1, \theta_2)\right)
\eeq
Now we observe that for the first two terms in the rhs of \eqref{four}
the action of the generator on the $\theta_1$ variable is the same as the action
of a reservoir generator on one site, with left parameter $\theta_L=0$ and right parameter $\theta_R=\theta_2$
(cf. \eqref{harmboundhid22}). For this generator, we know that the invariant measure is uniform
on $[0,\theta_2]$, which coincides with the conditional distribution $\Lambda^{2}_{0,1}(d\theta_1|\theta_2)$.
Therefore,
\beq
&&\int \Lambda^{2}_{0,1}(d\theta_1, d\theta_2)\int_0^1 \frac{du}{u}\left(f((1-u)\theta_1, \theta_2)-f(\theta_1, \theta_2)\right)
\nonumber\\
&+& \int \Lambda^{2}_{0,1}(d\theta_1, d\theta_2) \int_0^1 \frac{du}{u}\left(f((1-u)\theta_1+u\theta_2, \theta_2)-f(\theta_1, \theta_2)\right)
\nonumber\\
&=&
\int \Lambda^{2,2}_{0,1}(d\theta_2)\left(\int \Lambda^{2}_{0,1}(d\theta_1|\theta_2)\int_0^1 \frac{du}{u}\left(f((1-u)\theta_1, \theta_2)-f(\theta_1, \theta_2)\right)\right.
\nonumber\\
&+& \left.\int_0^1 \frac{du}{u}\left(f((1-u)\theta_1+u\theta_2, \theta_2)-f(\theta_1, \theta_2)\right)\right)=0
\eeq
Here we used the notation $\Lambda^{2,2}_{0,1}(d\theta_2)$ for the second marginal of the measure $\Lambda^{2}_{0,1}(d\theta_1, d\theta_2)$, and in the last step we used the invariance of the conditional distribution
$\Lambda^{2}_{0,1}(d\theta_1|\theta_2)$ for the reservoir generator with one site and left parameter zero, right parameter $\theta_2$.
Similarly, the action of the last two terms of the generator  in the rhs of \eqref{four} on the $\theta_2$ variable
is the same as the action of a reservoir generator on one site, with left parameter $\theta_L=\theta_1$ and right parameter $\theta_R=1$. As a consequence
\beq
&&\int \Lambda^{2}_{0,1}(d\theta_1, d\theta_2)\int_0^1 \frac{du}{u}\left(f(\theta_1, (1-u)\theta_2+u\theta_1)-f(\theta_1, \theta_2)\right)
\nonumber\\
&+&
\int\Lambda^{2}_{0,1}(d\theta_1, d\theta_2)\int_0^1 \frac{du}{u}\left(f(\theta_1, u+(1-u)\theta_2)-f(\theta_1, \theta_2)\right)
\nonumber\\
&=&
\int \Lambda^{2,1}_{0,1}(d\theta_1)\int \Lambda^{2}_{0,1}(d\theta_2|\theta_1)\left(\int_0^1 \frac{du}{u}\left(f(\theta_1, (1-u)\theta_2+u\theta_1)-f(\theta_1, \theta_2)\right)\right.
\nonumber\\
&+&
\left.\int_0^1 \frac{du}{u}\left(f(\theta_1, u+(1-u)\theta_2)-f(\theta_1, \theta_2)\right)\right)=0
\eeq
Here we used $\Lambda^{2,1}_{0,1}(d\theta_1)$ for the first marginal of the measure $\Lambda^{2}_{0,1}(d\theta_1, d\theta_2)$, and in the last step we used the invariance of the conditional distribution
$\Lambda^{2}_{0,1}(d\theta_2|\theta_1)$ for the reservoir generator with one site and left parameter $\theta_1$, right parameter $1$.
This finishes the proof of the case $N=2$.

The general case is now a straightforward generalization via the Markovian structure of the joint distribution of order statistics. Indeed, let
us call $\Lambda^{N}(d\theta_1, \ldots, d\theta_N)$ the joint distribution of the order statistics
of $N$ uniforms on $[\theta_L,\theta_R]$. Then conditional on $\theta_1, \ldots, \theta_{i-1}, \theta_{i+1}, \ldots, \theta_{N}$, the variable $\theta_i$ is uniform on
$[\theta_{i-1}, \theta_{i+1}]$ and
the terms in the generator acting on the variable $\theta_i$
exactly coincide with the action of the reservoir generator on a single
site with left reservoir parameter $\theta_{i-1}$ and right reservoir
parameter $\theta_{i+1}$, where we made the convention $\theta_0=\theta_L, \theta_{N+1}=\theta_R$.
\epr

The main reason why are able to identify the invariant measure for the hidden parameter model on the chain can be summarized as follows.
Consider the chain with left boundary parameter $\theta_L$ and right boundary parameter $\theta_R$.
Then the generator of the hidden parameter model has the form
\be\label{chaingengen}
\loc= \loc_1^{\theta_0} + \sum_{i=1}^N\loc_{i,i+1} + \loc_{N}^{\theta_{N+1}}
\ee
The action of this generator on the variable $\theta_i$  coincides with the action of
the generator \eqref{harmboundhid22} with $\theta_L=\theta_{i-1}, \theta_R=\theta_{i+1}$., 
In other words, the generator from \eqref{chaingengen} can be rewritten as
\be\label{chaingengengen}
\loc=\sum_{i=1}^N \loc^{\theta_{i-1},\theta_{i+1}}_i
\ee
where $\loc^{\theta_{i-1},\theta_{i+1}}_i$ is 
the generator \eqref{harmboundhid22} with $\theta_L=\theta_{i-1}, \theta_R=\theta_{i+1}$ acting on the $\theta_i$ variable.

Therefore,  if we call $\Lambda^1_{\theta_L,\theta_R}(d\theta)$ the invariant measure of the process with generator \eqref{harmboundhid22} describing the action of a left and right reservoir on a single site
(which is uniform for the generator \eqref{harmboundhid22}), then we can describe the invariant measure
of the generator \eqref{chaingengen} as follows. Let $\Lambda^{N}_{\theta_0,\theta_{N+1}}(d\theta_1,\ldots, d\theta_N)$ be a  probability measure  such
that its conditional distributions are given by
\[
\Lambda^{N}_{\theta_0,\theta_{N+1}}(d\theta_i|\theta_1,\ldots, \theta_{i-1}, \theta_{i+1}, \ldots \theta_{N})
= \Lambda^1_{\theta_{i-1},\theta_{i+1}}(d\theta_i)
\]
then $\Lambda^{N}_{\theta_0,\theta_{N+1}}(d\theta_1,\ldots, d\theta_N)$ is invariant for the generator \eqref{chaingengen}.

So we conclude that for every model for which one has this structure of the generator (i.e., \eqref{chaingengengen}), one can obtain its invariant measure once one has identified the invariant measure $\Lambda^1_{\theta_L,\theta_R}(d\theta)$ of the model with a single site between a left and right reservoir.
As we will see below, this is the case for the generalized harmonic models parametrized by $s>0$.

\subsection{The generalized harmonic model with parameter $s>0$: bulk generator}
The model from the previous section is
a special case of a one-parameter family of so-called generalized ``harmonic models''.
For these models, the measure
$M$ in \eqref{gengen} reads
\be\label{betaM}
M(du) =\frac{(1-u)^{2s-1}}{u} du
\ee
The corresponding discrete harmonic models  have the measure
\be\label{betaDisM}
\caM(k,n)= \frac{1}{k} \frac{\Gamma(n+1)\Gamma(n-k+2s)}{\Gamma(n+2s)\Gamma(n-k+1)} 
\ee
Considering the generalized  harmonic models
on a graph, we have then duality between the discrete model and the continuous model with duality functions
\be\label{dualfum}
D_m(\xi, x)= \prod_{i\in V} d_m(\xi_i, x_i)
\ee
with $d_m(k,x)=\frac{x^k\Gamma(2s)}{\Gamma(2s+k)}$
and self-duality of the discrete model with self-duality functions
\be\label{dualfuF}
D_F(\xi, \eta)= \prod_{i\in V} d_F(\xi_i, \eta_i)
\ee
with $d_F(k,n)= I(k\leq n) \frac{n!\Gamma(2s)}{(n-k)!\Gamma(2s+k)}$.
Moreover, the continuum model and the discrete model are Poisson intertwined.
See \cite{frass1}, \cite{frass2} for a proof of these dualities.
As a consequence, we have
the analogue of Theorem \ref{grandhalfdualitythm}, if one replaces
the duality functions $D_m$ and $D_F$ by
the ones in \eqref{dualfum}, resp. \eqref{dualfuF}.
The reversible product measures are now
given by products of Gamma distributions 
for the continuum model, and by products of discrete Gamma distributions for the discrete model.
As a consequence, we also have the analogue of Theorem \ref{mixthmharm} if one replaces the exponential $\nu_\theta$ marginals by the corresponding Gamma marginals for the continuum model, i.e.,
\be\label{gam2s}
\nu_\theta (dx)= \frac{x^{2s-1}}{\Gamma(2s)\theta^{2s}}e^{-x/\theta} dx
\ee
or by the corresponding discrete Gamma marginals for the discrete model
\[
\nu_\theta (n)
= 
\frac{1}{n!}\left(\frac{\theta}{1+\theta}\right)^n
\frac{\Gamma(2s+n)}{\Gamma(2s)}(1+\theta)^{-2s}
\]

\subsection{The generalized harmonic model with parameter $s>0$: boundary generator}
The boundary generator with reservoir parameter $\theta^*$ is \cite{frass2}
\beq\label{resers}
L_{\theta^*} f(x) &=&
\int_0^1 \frac{du}{u}(1-u)^{2s-1} \left(f((1-u)x)-f(x)\right)
\nonumber\\
&+&
\int_0^\infty  \frac{du}{u} e^{-u} \left(f(x+u\theta^*)-f(x)\right)
\eeq
The first term models the exit of mass to the reservoir, the second term models the input of mass from the reservoir.
Denoting by $\nu_{\theta}$ the Gamma distribution \eqref{gam2s}, the natural candidate intertwiner then reads
\be\label{colankogen}
\caG f(\theta)=
\int_0^\infty f(x) \nu_\theta(dx)
\ee
The candidate generator associated to a reservoir with parameter $\theta^*$ in the corresponding hidden parameter model
is  given by
\be\label{hiddensres}
\loc_{\theta^*} f(\theta)=
\int_0^1 \frac{du}{u}(1-u)^{2s-1}\left(f((1-u)\theta + u\theta^*)-f(\theta)\right)
\ee
We can now state the  analogue of the intertwining relation of Lemma \ref{interlem22}.
\bl
The boundary generator of the  generalized continuous harmonic process \eqref{resers}
and the boundary generator of the hidden parameter model \eqref{hiddensres}
are intertwined as
\be\label{intertws}
\caG L_{\theta^*} f= \loc_{\theta^*}\caG f
\ee
\el
\bpr
It suffices to see \eqref{intertws} for the functions $f(x)= x^n \frac{\Gamma(2s)}{\Gamma(2s+n)}$ (for all $n\in\N$).
We have
\begin{eqnarray}
&&
(\caG L_{\theta^*} f)(\theta)\nonumber\\
& = &
\int_0^\infty dx \frac{e^{-x/\theta}}{\theta} \left( \int_{0}^1 \frac{du}{u}(1-u)^{2s-1} \left(x^n(1-u)^n- x^n\right) +  \int_{0}^{\infty} \frac{du}{u} e^{-u} \left((x+u \theta^*)^n - x^n\right)\right)\frac{\Gamma(2s)}{\Gamma(2s+n)}
 \nonumber\\
 & = &
 \theta^n \int_{0}^1 \frac{du}{u} (1-u)^{2s-1}((1-u)^n-1) + \sum_{k=1}^n \theta^{n-k} (\theta^*)^k \frac{n!}{k!(n-k)!}\int_{0}^{\infty}\frac{du}{u} e^{-u}u^{k} \frac{\Gamma(2s+n-k)}{\Gamma(2s+n)}
  \nonumber\\
 & = &
 \theta^n \int_{0}^1 \frac{du}{u} (1-u)^{2s-1} ((1-u)^n-1) + \sum_{k=1}^n \theta^{n-k} (\theta^*)^k \frac{1}{k} \frac{\Gamma(2s+n-k)}{\Gamma(2s+n)} \frac{\Gamma(n+1)}{\Gamma(n+2s)}
 \end{eqnarray}
and we also have
 \begin{eqnarray}
 &&
(\loc_{\theta^*} \caG f)(\theta)\nonumber\\
& = &
\int_{0}^1 \frac{du}{u}(1-u)^{2s-1} \Big(  \int_{0}^{\infty} dx \frac{e^{-\frac{x}{(1-u) \theta+ u \theta^*}}}{(1-u) \theta+ u \theta^*}  x^n -  \int_{0}^{\infty} dx \frac{e^{-\frac{x}{\theta}}}{\theta}  x^n\Big)\frac{\Gamma(2s)}{\Gamma(2s+n)}
\nonumber\\
& = &
\int_{0}^1 \frac{du}{u}(1-u)^{2s-1} \Big(  ((1-u) \theta+ u \theta^*)^n -  \theta^n \Big)
\nonumber\\
 & = &
 \theta^n \int_{0}^1 \frac{du}{u} (1-u)^{2s-1} ((1-u)^n-1) + \sum_{k=1}^n \theta^{n-k} (\theta^*)^k \frac{1}{k} \frac{\Gamma(2s+n-k)}{\Gamma(2s+n)} \frac{\Gamma(n+1)}{\Gamma(n+2s)}
\end{eqnarray}
This completes the proof.
\epr

The generator of the hidden parameter model associated to two reservoirs acting on a single site is then given by

\be\label{hiddengeneral2ress}
\loc^{\theta_L,\theta_R} f(\theta)=
\loc_{\theta_L} f(\theta) +\loc_{\theta_R}f(\theta)
\ee
The following lemma identifies the invariant measure of 
$\loc^{0,1}$.
\bl\label{invmeasloc01}
The unique invariant measure of
the process with generator 
\eqref{hiddengeneral2ress} with
$\theta_L=0, \theta_R=1$ is equal to
the conditional distribution
\[
\Lambda^{0,1}(dy_1):=\nu_\theta\otimes \nu_\theta (dy_1|y_1+y_2=1)
\]
In particular, this is given by the Beta distribution
\[
Beta(2s, 2s)(d\theta)= \frac{\theta^{2s-1}(1-\theta)^{2s-1} }{B(2s,2s)} d\theta
\]
where $B(a,b)$ denotes the Beta function.
\el
\bpr
The action of $\loc^{0,1}$ on the
$\theta$ variable coincides with
the action of the generator
\eqref{gengen} with $M(du)= \frac{du}{u}(1-u)^{2s-1}$ when we start from
$y_1+y_2=1$, and consider the action on the $y_1$ variable.
The product measure $\nu_\theta\otimes \nu_\theta$ is a reversible measure, and
the event $y_1+y_2=1$ is invariant.
As a consequence, the
conditioned measure 
$\nu_\theta\otimes \nu_\theta(dy_1|y_1+y_2=1)$ is invariant.
Since $\nu_\theta$ is Gamma distributed
with shape parameter $2s$, then 
the conditional distribution
$\nu_\theta\otimes \nu_\theta(dy_1|y_1+y_2=1)$
is the symmetric Beta distribution with parameter $2s$.
\epr 
%
%
\subsection{Invariant measure of the generalized harmonic model on the chain}\label{invhidgen}
As a consequence of the intertwining of the boundary generator described in section 3.7, we obtain the analogue of Theorem \ref{contharminv} for the full set of boundary-driven generalized harmonic models with parameter $s$.
In order to understand the structure of the invariant measure in the
setting of the chain with left and right boundary reservoirs, we have to understand the invariant measure of the hidden parameter model.
We have seen in Lemma \ref{invmeasloc01} that the measure $\Lambda^{0,1}_1$ is
the   distribution $Beta(2s, 2s)[d\theta]$.
Let us call $B_{\theta_L,\theta_R}(d\theta)$ the corresponding recentered and rescaled distribution which is such that under this distribution 
\[
\frac{\theta-\theta_L}{\theta_R-\theta_L}
\]
is $Beta(2s,2s)$ distributed.

Then, following the line of argument of the proof of Theorem \ref{bofi} the invariant measure $\Lambda^N (d\theta_1, \ldots, d\theta_N)$ is such that its conditional
distributions are given by
\[
\Lambda^N (d\theta_i|\theta_0,\ldots, \theta_{i-1}, \ldots,\theta_{i+1},\ldots, \theta_{N+1})=
B_{\theta_{i-1}, \theta_{i+1}} (d\theta_i)
\]
This yields exactly the joint distribution
obtained in \cite{red}, i.e.,
\be
\label{invthetas}
\Lambda^N (d\theta_1, \ldots, d\theta_N)
=C(N,2s,\theta_L,\theta_R)\prod_{i=1}^{N+1}(\theta_i-\theta_{i-1})^{2s-1} \1(\theta_L\leq \theta_1\leq\ldots\leq \theta_N\leq \theta_R)
\ee
where $C(N,2s,\theta_L,\theta_R)$ is the normalization constant
$$
C(N,2s,\theta_L,\theta_R)
=
\frac{1}{(\theta_R-\theta_L)^{2s(N+1)-1}}
\frac{\Gamma(2s(N + 1))}{\Gamma(2s)^{N+1}}.
$$
We summarize the finding of this section in the
following theorem
\bt\label{bofi2}
The invariant measure of the hidden parameter model with generator
\begin{eqnarray}\label{ibolo}
&&
\caL f(\theta)\nonumber\\
&=& \int  \frac{du}{u}(1-u)^{2s-1} \Big(f(\theta- u\theta_1\delta_1 + u\theta_L \delta_1) - f(\theta)\Big)\nonumber\\
&+&
\sum_{i=1}^N 
\int_0^1 \frac{du}{u}(1-u)^{2s-1} \left(f(\theta- u\theta_i\delta_i + u\theta_{i+1} \delta_i)+f(\theta- u\theta_{i+1}\delta_{i+1} + u\theta_i \delta_{i+1})-2f(\theta)\right)
\nonumber\\
&+& \int  \frac{du}{u}(1-u)^{2s-1} \Big(f(\theta- u\theta_N\delta_N + u\theta_R \delta_N) - f(\theta)\Big)
\end{eqnarray}
is the measure \eqref{invthetas}.
As a consequence, the invariant measure of the boundary driven continuous harmonic model with parameter $s$ defined by the generator 
\begin{eqnarray}\label{iboli}
&&
L f(\zeta)\nonumber\\
&=& \int_0^1 \frac{du}{u}(1-u)^{2s-1} \left(f(\zeta -u\delta_1)-f(\zeta)\right)
+
\int_0^\infty  \frac{du}{u} e^{-u} \left(f(\zeta+u\theta_L \delta_1)-f(\zeta)\right)\nonumber\\
&+&
\sum_{i=1}^N 
\int_0^1 \frac{du}{u}(1-u)^{2s-1} \left(f(\zeta- u\zeta_i\delta_i + u\zeta_{i} \delta_{i+1})+f(\zeta- u\zeta_{i+1}\delta_{i+1} + u\zeta_{i+1} \delta_{i})-2f(\zeta)\right)
\nonumber\\
&+& \int_0^1 \frac{du}{u}(1-u)^{2s-1} \left(f(\zeta -u\delta_N)-f(\zeta)\right)
+
\int_0^\infty  \frac{du}{u} e^{-u} \left(f(\zeta+u\theta_R \delta_N)-f(\zeta)\right)
\end{eqnarray}
is a mixture of product of Gamma distributions with mixing measure \eqref{invthetas}.
\et
\subsection{General redistribution rules and reservoirs}
We close this section by investigating intertwining for the general mass redistribution model with generator
\eqref{gengen}. We also discuss a general definition of reservoirs which is naturally associated to the the general mass redistribution model.

We assume  that the measure $M$ in \eqref{gengen} is chosen in such a way that the corresponding process has a one-parameter family of reversible product measures
with marginals denoted by $\nu_\theta(dx)$. E.g. for the choice
$M(du)=(1/u) du$, $\nu_\theta(dx)= \frac{1}{\theta}e^{-x/\theta} dx$; for the choice
$M(du)=(1/u)(1-u)^{2s-1} du$, $\nu_\theta(dx)= \frac{1}{\theta}e^{-x/\theta} x^{2s-1} dx$.
The natural boundary generator with reservoir parameter $\theta^*$ is given by
\beq\label{resergen}
L_{\theta^*} f(x) &=&
\int_0^1 M(du)
 \left(f((1-u)x)-f(x)\right)
\nonumber\\
&+&
\int_0^{\infty} \nu_{\theta^*}(dy)\int_0^1  M(du)\left(f(x+u y)-f(x)\right)
\eeq
As in the KMP process, this choice of the reservoir
in inspired by the idea that a site interact with the reservoir as it does with the bulk sites it is connected to; however the energy of the reservoirs
is random and sampled from the distribution $\nu_{\theta^*}$.

In the corresponding hidden parameter model, the candidate generator associated to a reservoir with parameter $\theta^*$ is then given by
\be\label{hiddengeneralres}
\loc_{\theta^*} f(\theta)=
\int_0^1 M(du)\left(f((1-u)\theta + u\theta^*)-f(\theta)\right)
\ee
and the generator of the hidden parameter model associated to two reservoirs acting on a single site is 
\be\label{hiddengeneral2res}
\loc^{\theta_L,\theta_R} f(\theta)=
\loc_{\theta_L} f(\theta) +\loc_{\theta_R}f(\theta)
\ee
The following lemma identifies the invariant measure of 
$\loc^{0,1}$ in terms of the measure $\nu_{\theta}$.
\bl\label{invmeasloc01}
The unique invariant measure of
the process with generator 
\eqref{hiddengeneral2res} with
$\theta_L=0, \theta_R=1$ is equal to
the conditional distribution
\[
\Lambda^{0,1}(dx_1):=\nu_\theta\otimes \nu_\theta (dx_1|x_1+x_2=1)
\]
In particular, the latter does not depend on $\theta$.
\el
\bpr
The action of $\loc^{0,1}$ on the
$\theta$ variable coincides with
the action of the generator
\eqref{gengen} when we start from
$y_1+y_2=1$, and consider the action on the $y_1$ variable.
By assumption, $\nu_\theta\otimes \nu_\theta$ is a reversible measure, and
the event $y_1+y_2=1$ is invariant.
As a consequence, the
conditioned measure 
$\nu_\theta\otimes \nu_\theta(dy_1|y_1+y_2=1)$ is invariant.
\epr 

We introduce the natural candidate intertwiner as
\be\label{colankogen}
\caG f(\theta)=
\int_0^\infty f(x) \nu_\theta(dx)
\ee
and its tensorization
\be\label{tensor}
\caG f((\theta_i)_{i\in V})
=
\int 
 f((x_i)_{i\in V}) \bigotimes\nu_{\theta_i}(\prod_i dx_i)
\ee
To discuss intertwining for the boundary-driven model
we need to establish conditions guaranteeing that
\be\label{intertwgeneral}
\caG L_{\theta^*}= \loc_{\theta^*}\caG
\ee
In order to obtain the intertwining \eqref{intertwgeneral} we make the following natural scaling assumption on the measure $\nu_\theta$:
\be\label{nuthetascaling}
\int \nu_\theta (dx) x^n= R_n \theta^n
\ee
Here we implicitly assumed that all the moments are finite. We moreover assumed that the measures $\nu_\theta$ are uniquely determined by their moments.
Then we have the following.
\bl\label{intertwlemgen}
The intertwining relation \eqref{intertwgeneral} is satisfied if and only if for all $n$ and $k\in\{0,\ldots, n\}$ we have
\be\label{momrec}
R_{n-k}R_k \int_0^1 u^k \ M(du) 
= R_n \int_0^1 M(du) u^k (1-u)^{n-k}
\ee
\el
\bpr
We start from \eqref{intertwgeneral} and
fill in  the function $f(x)= x^n$.
Then the lhs equals
\beq
(\caG (L_{\theta^*} f))(\theta)
&=&
R_n \theta^n \int_0^1 M(du) ((1-u)^n-1)
\nonumber\\
&+&
\sum_{k=1}^n \theta^{n-k}(\theta^*)^k R_{n-k}R_k {n\choose k} \int_0^1 M(du) u^k
\eeq
The rhs equals
\beq
(\loc_{\theta^*}(\caG f)) (\theta)
&=&
R_n \theta^n \int_0^1 ((1-u)^n-1) M(du) 
\nonumber\\
&+&
\sum_{k=1}^n \theta^{n-k}(\theta^*)^k R_n
{n\choose k}\int_0^1 u^k (1-u)^{n-k} M(du) 
\eeq 
Because both expressions have to be equal for all values of $\theta, \theta^*$, we obtain
\eqref{momrec}
\epr 

The following corollary then proves that
\eqref{momrec} is satisfied for the harmonic model with $s>0$.
\bc
Let $\nu_\theta$ be the Gamma distribution
with scale parameter $\theta$ and shape parameter $2s$, i.e., the probability measure given in \eqref{gam2s}, and assume
\be\label{Mdus}
M(du)= \frac{(1-u)^{2s-1}}{u} du
\ee
Then the moment relation
\eqref{momrec} is satisfied.
\ec 
\bpr
For the measure \eqref{Mdus} we have $R_n = \frac{\Gamma(2s+n)}{\Gamma(2s)}$. Therefore,
using \eqref{Mdus},
the lhs of \eqref{momrec} equals
\[
R_{n-k}R_k \int_0^1 u^k \ M(du) 
= \frac{\Gamma (n-k +2s)}{\Gamma (2s)} \frac{\Gamma(k+2s)}{\Gamma(2s)} \frac{\Gamma(k)\Gamma(2s)}
{\Gamma(k+2s)}
\]
and the rhs of \eqref{momrec} equals
\[
R_n \int_0^1 M(du) u^k (1-u)^{n-k}=\frac{\Gamma(n+2s)}{\Gamma(2s)}\frac{\Gamma(k)
\Gamma(n-k+2s)}{\Gamma(n+2s)}
\]
Hence both expressions are indeed equal and the relation
\eqref{momrec} is satisfied.
\epr 
\br
As a consequence of the above corollary we deduce that the boundary-driven generalized harmonic models with reservoirs
\eqref{resers} and the boundary-driven generalized harmonic models with reservoirs 
\beq\label{resfin}
L_{\theta^*} f(x) &=&
\int_0^1 \frac{du}{u}(1-u)^{2s-1}
 \left(f((1-u)x)-f(x)\right)
\nonumber\\
&+&
\int_0^{\infty} \nu_{\theta^*}(dy)\int_0^1  \frac{du}{u}(1-u)^{2s-1}\left(f(x+u y)-f(x)\right)
\eeq
have in turn the same hidden parameter model. Therefore they will also have the {\em same} stationary measure, which in the case of the chain can be explitely characterized as
a mixture of product of Gamma distributions with mixing measure \eqref{invthetas}.
\er

\section{Poisson intertwining}
\label{section4}
The class of models for which the non-equilibrium steady state is a mixture of product measures is not limited to models of KMP or harmonic type. In this section we consider the boundary driven symmetric inclusion process (SIP) and we prove that the stationary measure is a mixture of Poisson product measures. This is a different situation compared to the previous sections, because the Poisson product measures are not the stationary measures of the SIP. The stationary measures of SIP are product of discrete Gamma distributions, which are however themselves mixtures of Poisson product measures.

To prove the result for SIP,  we use a simple Poissonian intertwiner of the classical creation and annihilation operators which transforms the boundary generators into the boundary generators of the Brownian energy process (BEP).
Using this same intertwiner, we also revisit the simplest example of independent random walkers, by which we then recover the propagation of Poisson product measures, which is a version of Doob's theorem \cite{dmpres}, or alternatively,  of the random displacement theorem in point process theory \cite{last}.

\subsection{Boundary driven SIP}
First, the SIP on two sites is the Markov process on $\N^2$ with generator
\be\label{twosip}
L f(n_1,n_2)= n_1(2s+ n_2) (f(n_1-1,n_2+1)-f(n_1,n_2)) + n_2(2s+ n_1) (f(n_1+1,n_2-1)-f(n_1,n_2))
\ee
Given a vertex set $V$ and irreducible edge weights $p(i,j)$, we define the SIP as the Markov process
with generator
\be\label{sipgen}
\sum_{ij} p(i,j) L_{ij}
\ee
where, as usual, $L_{ij}$ is the generator \eqref{twosip} acting on the variables $\eta_i,\eta_j$.
The boundary generator is given by
\be\label{boundgensip}
L^{\alpha,\gamma} f(n)= \alpha (2s+n) (f(n+1)-f(n)) + \gamma n (f(n-1)-f(n))
\ee
we assume $\alpha<\gamma$, in that case $L^{\alpha, \gamma}$ admits a unique stationary measure which is the discrete Gamma distribution \eqref{discgam} with $\theta= \frac{\alpha}{\gamma-\alpha}$.

The boundary driven model with boundary reservoirs  is then given by
\be\label{boundsip}
\sum_{ij} p(i,j) L_{ij} + \sum_{i\in V} c(i) L^{\alpha_i,\gamma_i}_i
\ee
where $L_i$ denotes the generator \eqref{boundgensip} acting on the variable $\eta_i$.

We first rewrite the boundary generators in terms of creation and annihilation operators.
The latter are defined as acting on a function $f:\N\to\R$ via
\beq\label{anni}
a f(n) &=& n f(n-1)
\nonumber\\
a^\dagger f(n) &=& f(n+1)
\eeq
where in \eqref{anni} it is understood $af(0)=0$.
We denote by $a_i$, resp. $a^\dagger_i$ these operators acting on the variable $\eta_i$. Then these
operators satisfy the conjugate Heisenberg algebra commutation relations, i.e.
\beq\label{commu}
[a_i,a_j]=[a^\dagger_i, a^\dagger_j]=0, \ [a_i^\dagger, a_j] = \delta_{i,j}
\eeq
and we can rewrite the boundary generator \eqref{boundgensip} as
\be\label{boundcre}
L^{\alpha,\gamma} = 2s\alpha (a^\dagger -I) +\alpha (aa^\dagger a^\dagger - aa^\dagger) + \gamma (a- aa^\dagger)
\ee
We  first define the Poisson intertwining which turns the operators $a,a^\dagger$ into differential operators.
\bl\label{poissin}
Define, for $f:\N\to\R$
and $z\geq 0$
\be\label{poissindef}
\caG f(z)=\sum_{n=0}^\infty \frac{z^n}{n!} e^{-z} f(n)
\ee
Then we have
\beq\label{interpois}
\caG a f &=& A\caG f, \ \caG a^\dagger f = A^\dagger \caG f
\eeq
with
\beq\label{AA}
Af(z) = z f(z)
,\
A^\dagger f(z)= f'(z) + f(z)
\eeq
As a consequence
\be\label{boundgenA}
\caG L^{\alpha,\gamma} =  \loc^{\alpha,\gamma} \caG
\ee
with
\be\label{difboundgen}
\loc^{\alpha,\gamma} = \alpha z \partial_z^2 + (2s\alpha - (\gamma-\alpha)z) \partial_z
\ee
\el
\bpr
The intertwinings \eqref{interpois} follow from a direct computation.
Then \eqref{boundgenA} follows from \eqref{boundcre} and \eqref{interpois}.
\epr
Notice that $\caG f(z)=\int f(n) \pi_z(dn)$ where $\pi_z$ is the Poisson measure with parameter $z$.
We extend as usual the intertwining $\Lambda$ by tensorization, i.e., for $f:\N^V\to\R$
\[
\caG f(\zeta)=
\int f(\eta)( \otimes \pi_{\zeta_i}) (d\eta)
\]
then we have the following intertwining result.
\bt\label{sipinter}
The boundary driven SIP with generator \eqref{boundsip} is Poisson intertwined with the boundary driven BEP process with generator
\be\label{boundbep}
\loc= \sum_{i,j\in V} p(i,j) \loc_{ij} + \sum_{i\in V} c(i)\loc^{\alpha_i,\gamma_i}_i
\ee
Here $\loc_{ij}$ is the single edge generator of the Brownian enery process (BEP), given by
\be\label{bepsigen}
\loc_{ij}= \zeta_i\zeta_j(\partial_{i}-\partial_j)^2- 2s(\zeta_i-\zeta_j) (\partial_i-\partial_j)
\ee
Here $\partial_i$ denotes partial derivative w.r.t.\ $\zeta_i$,
and $\loc^{\alpha_i,\gamma_i}_i$ is \eqref{difboundgen} acting on the variable $\zeta_i$.

As a consequence we have the following result on propagation of Poisson product measures. If we start the boundary driven SIP
from the product Poisson measure $\otimes_{i\in V} \pi_{\zeta_i}$ then we have
\be
(\otimes_{i\in V} \pi_{\zeta_i}) S(t)=\E_{\zeta} \left(\otimes_{i\in V} \pi_{\zeta_i(t)}\right)
\ee
where $\zeta(t)$ evolves according to the generator $\loc$ in \eqref{boundbep}.
As a further consequence the unique stationary measure of the boundary driven SIP is a mixture of Poisson measures with
\[
\int(\otimes_{i\in V} \pi_{\zeta_i})\Xi(\prod_{i\in V}d\zeta_i)
\]
where the mixture measure $\Xi$ is the unique stationary measure of the process with generator
\eqref{boundbep}.
\et
\bpr
The intertwining follows from the combination of Lemma \ref{poissin}  with the fact that
the single edge generators of SIP and BEP are Poisson intertwined see e.g. \cite{fede}, or \cite{book}, i.e., for all $ij$
\[
\caG L_{ij}= \loc_{ij} \caG
\]
\epr
\subsection{Independent random walkers}
The independent random walk process on two sites is given by the generator
\[
L_{12}f(n_1,n_2)= n_1 (f(n_1-1,n_2+1)-f(n_1,n_2)) + n_2(f(n_1-1,n_2+1)-f(n_1,n_2))
\]
and the boundary generator
\[
L f(n)= \alpha (f(n+1)-f(n)) + \gamma n (f(n-1)-f(n))
\]
Then the full boundary driven model for independent random walkers reads as follows.
\be\label{boundirw}
L= \sum_{i,j\in V}p(i,j) L_{ij}  + \sum_{i\in V} c(i) L^{\alpha_i,\gamma_i}_i
\ee
In terms of the creation and annihilation operators the generators read
\be\label{irwsingle}
L_{ij}= -(a_i-a_j)(a_i^\dagger- a_j^\dagger)
\ee
for the single edge generator
and
\be\label{irwbound}
L^{\alpha_i,\gamma_i}_i= \alpha_i (a_i^\dagger-I) + \gamma_i (a_i-a_ia_i^\dagger)
\ee
for the boundary generator.
Therefore, using Lemma \ref{poissin}, we obtain that the boundary generator is Poisson intertwined with the operator
\be\label{bopopo}
\loc^{\alpha_i,\gamma_i}_i= \alpha_i (A_i^\dagger-I) + \gamma_i (A_i-A_iA_i^\dagger)=
(\alpha_i-\gamma_i z_i)\partial_i
\ee
and the single edge generator is intertwined with the operator
\be\label{singleirw}
\loc_{ij}= -(\zeta_i-\zeta_j)(\partial_i-\partial_j)
\ee
Notice that $\loc^{\alpha_i,\gamma_i}_i$ and $\loc_{ij}$ are first order differential operators and therefore the process build from them is a deterministic system of ODEs.
We then immediately obtain the following analogue of Theorem \ref{sipinter}
\bt\label{irwinter}
The boundary driven independent random walkers with generator \eqref{boundirw} is Poisson intertwined with the boundary driven deterministic process with generator
\be\label{boundbep22}
\loc= \sum_{i,j\in V} p(i,j) \loc_{ij} + \sum_{i\in V}\loc^{\alpha_i,\gamma_i}_i
\ee
Here $\loc_{ij}$ is the single edge generator \eqref{singleirw}
and $\loc^{\alpha_i,\gamma_i}_i$ is \eqref{bopopo}.

As a consequence we have the following. When we start the boundary driven SIP
from the product Poisson measure $\otimes_{i\in V} \pi_{\zeta_i}$ then we have
\be
(\otimes_{i\in V} \pi_{\zeta_i}) S(t)= \otimes_{i\in V} \pi_{Z^{\zeta}_i(t)}
\ee
where $Z^\zeta(t)$ evolves according to the deterministic generator $\loc$ in \eqref{boundbep22}.
As a further consequence the unique stationary measure of the boundary driven independent random walkers is a Poisson product measure
\[
\otimes_{i\in V} \pi_{\zeta^*_i}
\]
$\zeta^*$ is the unique fixed point of the deterministic system $Z^\zeta(t)$.
\et

\bibliographystyle{plain}
\bibliography{refsinter}

\end{document}